\documentclass{article}

\usepackage[english]{babel}
\usepackage[utf8x]{inputenc}
\usepackage[numbers]{natbib}
\usepackage{authblk}
\usepackage{url}

\usepackage{mathtools}
\usepackage{amsmath}
\usepackage{amssymb}
\usepackage{amsthm}
\usepackage{mathabx}

\usepackage{enumerate}	
\usepackage{xcolor}

\usepackage{multicol}
\usepackage{multirow}
\usepackage{makecell}
\usepackage{booktabs}
\usepackage{arydshln}
\usepackage{caption}
\usepackage{subcaption}
\usepackage{varwidth}
\usepackage{algorithm}
\usepackage{algpseudocode}

\usepackage{tikz}
\usepackage{tikz-3dplot}
\usepackage{pgf}
\usepackage{pgfplots}
\pgfplotsset{compat=1.17}
\usetikzlibrary{arrows.meta}
\usetikzlibrary{positioning}
\usetikzlibrary{calc}

\usepackage[capitalize,noabbrev,nameinlink]{cleveref}

\renewcommand{\arraystretch}{1.5}
\DeclareMathOperator{\argmin}{argmin}

\DeclareMathOperator{\epi}{epi}
\DeclareMathOperator{\spn}{span}
\DeclareMathOperator{\aff}{aff}
\DeclareMathOperator{\ri}{ri}

\algrenewcommand\algorithmicrequire{\textbf{Input:}}
\algrenewcommand\algorithmicensure{\textbf{Output:}}

\newtheorem{theorem}{Theorem}
\newtheorem{proposition}[theorem]{Proposition}
\newtheorem{lemma}[theorem]{Lemma}
\newtheorem{corollary}[theorem]{Corollary}

\theoremstyle{definition}
\newtheorem{definition}[theorem]{Definition}
\newtheorem{example}{Example}

\crefformat{program}{#2Program~(#1)#3}
\Crefformat{program}{#2Program~(#1)#3}
\crefmultiformat{program}{#2Programs~(#1)#3}{ and~#2(#1)#3}{, #2(#1)#3}{, and~#2(#1)#3}
\Crefmultiformat{program}{#2Programs~(#1)#3}{ and~#2(#1)#3}{, #2(#1)#3}{, and~#2(#1)#3}

\crefformat{progeq}{#2(#1)#3}
\crefmultiformat{progeq}{#2(#1)#3}{ and~#2(#1)#3}{, #2(#1)#3}{, and~#2(#1)#3}

\newcommand\numberthis{\addtocounter{equation}{1}\tag{\theequation}}

\newcommand{\etal}{et al.}
\newcommand{\ie}{\emph{i.e.}~}
\newcommand{\eg}{\emph{e.g.}~}
\newcommand{\st}{\text{s.t.~}}

\newcommand{\tr}[1]{{#1}^\top}

\newcommand{\setV}{\mathcal{V}}
\newcommand{\setA}{\mathcal{A}}
\newcommand{\setT}{\mathcal{T}}
\newcommand{\setK}{\mathcal{K}}
\newcommand{\setX}{\mathcal{X}}

\newcommand{\setP}{\mathcal{P}}
\newcommand{\setQ}{\mathcal{Q}}
\newcommand{\setreal}{\mathbb{R}}

\newcommand{\funcT}{\mathbf{T}}
\newcommand{\funcW}{\mathbf{W}}
\newcommand{\funcR}{\mathbf{R}}

\newcommand{\funcG}{\mathbf{G}}

\newcommand{\sumk}{\sum_{k \in \setK}}
\newcommand{\sumpk}{\sum_{\pk \in \setP^k}}

\newcommand{\cx}{\tr{c}x}
\newcommand{\tx}{\tr{t}x}
\newcommand{\ctx}{\tr{(c + t)}x}
\newcommand{\by}{\tr{b}y}
\newcommand{\Ax}{Ax}
\newcommand{\Ay}{\tr{A}y}
\newcommand{\xt}{\tr{x}t}
\newcommand{\xA}{x_{\setA_1}}
\newcommand{\xB}{x_{\setA_2}}
\newcommand{\txA}{\tr{t}\xA}

\newcommand{\ok}{{o^k}}
\newcommand{\dk}{{d^k}}
\newcommand{\bk}{{b^k}}
\newcommand{\xk}{x^k}
\newcommand{\yk}{y^k}
\newcommand{\cxk}{\tr{c}\xk}
\newcommand{\txk}{\tr{t}\xk}
\newcommand{\ctxk}{\tr{(c + t)}\xk}
\newcommand{\byk}{\tr{(\bk)}\yk}
\newcommand{\Axk}{A\xk}
\newcommand{\Ayk}{\tr{A}\yk}
\newcommand{\xkA}{\xk_{\setA_1}}
\newcommand{\xkB}{\xk_{\setA_2}}
\newcommand{\txkA}{\txk_{\setA_1}}

\newcommand{\tw}{\tr{t}w}
\newcommand{\wt}{\tr{w}t}

\newcommand{\wk}{w^k}
\newcommand{\twk}{\tr{t}\wk}

\newcommand{\hx}{\hat{x}}
\newcommand{\hxA}{\hx_{\setA_1}}
\newcommand{\hw}{\hat{w}}
\newcommand{\hwk}{\hw^k}
\newcommand{\cw}{\check{w}}
\newcommand{\cwk}{\cw^k}

\newcommand{\pk}{{p^k}}
\newcommand{\zk}{z^k}
\newcommand{\zkp}{\zk_\pk}

\newcommand{\xp}{\hx^p}
\newcommand{\xq}{\hx^q}

\newcommand{\xpk}{\hx^\pk}
\newcommand{\xpka}{\xpk_a}
\newcommand{\cpk}{\tr{c}\xpk}
\newcommand{\tpk}{\tr{t}\xpk}
\newcommand{\ctpk}{\tr{(c + t)}\xpk}
\newcommand{\xpkA}{\xpk_{\setA_1}}

\title{Asymmetry in the Complexity of the Multi-Commodity Network Pricing Problem}
\author[1]{Quang Minh Bui}
\author[1]{Margarida Carvalho}
\author[2]{Jos\'e Neto}
\affil[1]{CIRRELT and D\'epartement d'informatique et de recherche op\'erationnelle, Universit\'e de Montr\'eal}
\affil[2]{T\'el\'ecom SudParis, Institut Polytechnique de Paris}
\date{}

\begin{document}

\maketitle

\begin{abstract}
The network pricing problem (NPP) is a bilevel problem, where the leader optimizes its revenue by deciding
on the prices of certain arcs in a graph, while expecting the followers (also known as the commodities)
to choose a shortest path based on those prices.
In this paper, we investigate the complexity of the NPP with respect to two parameters:
the number of tolled arcs, and the number of commodities. We devise a simple algorithm showing that
if the number of tolled arcs is fixed, then the problem can be solved in polynomial time with respect to
the number of commodities. In contrast, even if there is only one commodity, once the number of tolled arcs is not fixed,
the problem becomes NP-hard. We characterize this asymmetry in the complexity with a novel property named strong bilevel feasibility.
Finally, we describe an algorithm to generate valid inequalities to the NPP based on this property, accommodated with
numerical results to demonstrate its effectiveness in solving the NPP with a high number of commodities.
\end{abstract}

\setlength{\textfloatsep}{10pt plus 1.0pt minus 2.0pt}
\section{Introduction}
In the most basic version of the network pricing problem (NPP), we have an optimization problem involving a graph and two decision makers:
the leader and the follower. First, the leader adjusts the prices of some arcs in the graph (called tolled arcs).
Subsequently, the follower finds the shortest path between its origin and destination according to the prices set by the leader.
The objective of the leader is to maximize the overall revenue which depends on whether the follower uses the tolled arcs or not.
To accomplish this goal, the leader wants to set the prices as high as possible, but not too high so that the follower is still incentivized to use
the leader's service. In this work, we are interested in the multi-commodity variant of the NPP, in which there
are multiple followers, each one has its own origin/destination, and each one chooses its own shortest path across the network.

\paragraph{Motivation}
The single-commodity NPP was first introduced by Labb\'e~\etal~\cite{labbe1998}. Since the follower's problem is the shortest path problem,
it can be written as a linear program. Using the Karush-Kuhn-Tucker (KKT) conditions, the NPP and its multi-commodity variant
can be formalized as a mixed-integer linear program (MILP) \citep{brotcorne2001}. Bui~\etal~\cite{bui2022} summarized and extended various techniques
to derive an MILP for the NPP, which covers path-based formulations \citep{didibiha2006, bouhtou2007}
and preprocessing \citep{bouhtou2007, vanhoesel2003}. Non-MILP methods include multipath enumeration \citep{brotcorne2011} and tabu search \citep{brotcorne2012}.

The single-commodity NPP has been proven to be NP-hard \citep{roch2005}. The complexity of the NPP is tied to the number of tolled arcs in the graph.
Indeed, if we consider the case where there is only one tolled arc and multiple commodities, then the NPP can be solved in polynomial time \citep{labbe2021}.
Thus, an asymmetry exists in the complexity of the NPP between the number of tolled arcs (denoted as \(|\setA_1|\)) and the number of commodities (denoted as \(|\setK|\)).
For the case with multiple commodities and multiple tolled arcs, van Hoesel~\etal~\cite{vanhoesel2003} proved that if \(|\setA_1|\) is fixed, then the NPP
can be solved in polynomial time with respect to \(|\setK|\). Specifically, the NPP can be decomposed into \(|\setK|^{f(|\setA_1|)} g(|\setA_1|)\) linear programs
of size \(h(|\setA_1|)\) each, where \(f(\cdot), g(\cdot), h(\cdot)\) are functions of exponential order.
Even with very small \(|\setA_1|\), the number of linear programs required to be solved is enormous, thus the algorithm described in \citep{vanhoesel2003} has no practical use.

\paragraph{Contributions and Paper Organization}
In this paper, we improve the result regarding the asymmetry in the complexity with the following theorem:
\begin{theorem}\label{thm:multi-comm-poly}
	If the number of tolled arcs \(|\setA_1|\) is fixed, then the multi-commodity NPP can be solved in polynomial time with respect to
	the number of commodities \(|\setK|\), specifically by solving \((|\setK| + 1)^{|\setA_1|}\) linear programs, each with size polynomial in \(|\setA_1|\).
\end{theorem}
Compared to \citep{vanhoesel2003}, the number of linear programs is much smaller.
Proving \cref{thm:multi-comm-poly} is the main focus of \cref{sec:asymmetry}.
In that section, we first revise the complexity of the single-commodity and the single-tolled-arc cases.
Then, we extend the asymmetry in the complexity to the multi-commodity case in an intuitive way by using \emph{reaction plots}.
Finally, we provide a rigorous proof of this asymmetry via a new reformulation of the NPP called the \emph{conjugate model}.

With the key result proven, we aim to exploit this asymmetry in a practical manner.
In \cref{sec:strongbf}, we introduce and utilize a new concept named \emph{strong bilevel feasibility} to generate cuts to existing MILP formulations of the NPP.
Strong bilevel feasibility is a property of a composition of paths across multiple followers. We show that there is always a solution of the NPP
that is strongly bilevel feasible, thus it is sufficient to enumerate only the strongly bilevel feasible points.
This is analogous to the fact that it is sufficient to enumerate only the extreme points of the feasible set to solve a linear program.
We also derive the characteristics of strong bilevel feasibility with the help of convex conjugates.

\Cref{sec:experiments} describes the cut generation procedure using strong bilevel feasibility.
We conduct numerical experiments to demonstrate that the asymmetry in the complexity is relevant in practice and that these cuts
are effective in accelerating the solution of the problem instances with a very high number of commodities.
\Cref{sec:conclusion} concludes the paper.

\section{Asymmetry in the Complexity}\label{sec:asymmetry}
\Cref{thm:multi-comm-poly} implies that the multi-commodity NPP scales differently with respect to two parameters:
the number of commodities or the number of tolled arcs.
To understand the difference between these two parameters, it is helpful to build some intuition from simpler cases:
the single-commodity case (\cref{ssec:single-comm}), and the single-tolled-arc case (\cref{ssec:single-toll}).
Then, in \cref{ssec:reaction}, we introduce an
illustration tool called reaction plot, and expand the intuition to the multi-commodity variant by superimposing these reaction plots
on top of each other. Finally, we materialize this intuition and provide a proof of \cref{thm:multi-comm-poly} in \cref{ssec:conjugate}.

\subsection{Single-Commodity Case}\label{ssec:single-comm}
First, we need to describe the NPP in its single-commodity form.
Let us consider a directed graph \(G = (\setV, \setA)\) where \(\setV\) and \(\setA\) are the set of vertices and the set of arcs, respectively.
The leader controls the costs of a subset of arcs \(\setA_1 \subsetneq \setA\), designated as tolled arcs.
All the other arcs form the set of toll-free arcs \(\setA_2 = \setA \setminus \setA_1\). The cost of a toll-free arc \(a \in \setA_2\) is always \(c_a \geq 0\),
while the cost of a tolled arc \(a \in \setA_1\) is \(c_a + t_a\) where \(t_a \geq 0\) is the price set by the leader (we only consider non-negative toll prices).

The follower wants to travel from an origin node \(o \in \setV\) to a destination \(d \in \setV\) via the shortest path across the graph.
After the follower has chosen the path, the leader collects a revenue equal to the sum of \(t_a\) over the arcs that the follower uses.
If there are multiple shortest paths, the follower chooses the path that produces the highest revenue for the leader (optimistic assumption).
To prevent the leader to extract infinite revenue from the follower, we assume that there always exists a toll-free path from \(o\) to \(d\).

Let \(x \in \{0, 1\}^{\setA}\) represent the selection of arcs corresponding to the shortest path chosen by the follower, where \(x_a = 1\) if arc \(a \in \setA\) is in the path
and \(x_a = 0\) otherwise. The NPP can be formulated as a bilevel program:
\begin{equation*}
	\max_{t,x} \lbrace \tx \mid  t \in \setT,  x \in \funcR(t) \rbrace
\end{equation*}
where \(\funcR(t)\) is the reaction set of the follower, defined by:
\begin{equation*}
	\funcR(t) =  \argmin_x \lbrace \ctx \mid \Ax = b,   x \geq 0\rbrace.
\end{equation*}
The set \(\setT = \{t \geq 0 \mid t_a = 0, \forall a \in \setA_2\}\) contains all vectors of feasible toll prices. The matrix \(A\) is
the incidence matrix of the graph, while \(b\) is the source-sink vector of the shortest path problem: \(b_o = 1\), \(b_d = -1\),
and \(b_i = 0\) for all other nodes \(i \in \setV\).

The canonical way to solve the NPP is to convert the bilevel program into a single-level reformulation
using the KKT conditions (with strong duality) \citep{brotcorne2001}:
\begin{subequations}
	\label[program]{prog:1comm-1level}
	\begin{align}
		\max_{t,x,y}\  & \tx             \label[progeq]{eq:1comm-1level-obj} \\
		\st            & t \in \setT,                                        \\
		               & \Ax = b,                                            \\
		               & \Ay \leq c + t,                                     \\
		               & \ctx = \by,     \label[progeq]{eq:1comm-1level-sd}  \\
		               & x \geq 0. 		 \label[progeq]{eq:1comm-1level-end}
	\end{align}
\end{subequations}
To solve this problem as an MILP, we need to linearize the bilinear term \(\tx\) appearing in the objective function \cref{eq:1comm-1level-obj}
and the strong duality constraint \cref{eq:1comm-1level-sd}. Since the optimal solution of the shortest path problem needs \(x\) to be binary, the
McCormick envelope \citep{mccormick1976} can be applied by replacing \(t_a x_a\) with \(s_a\) and appending for all $a \in \setA_1$:
\begin{align*}
	0 \leq s_a \leq M_a x_a, \qquad 0 \leq t_a - s_a \leq M_a (1 - x_a), \qquad  x_a \in \{0, 1\},
\end{align*}
to \cref{prog:1comm-1level} \citep{brotcorne2001}. Appropriate values of the big-M parameters \(M_a\) can be found in Dewez~\etal~\cite{dewez2008}.

In this form, we can see that the variables \(x_a,\ a \in \setA_1\), are the only binary variables in the MILP.
Thus, there is an indication that the difficulty of the NPP relies on the number of tolled arcs~\(|\setA_1|\).
We denote this truncated vector \(\xA \in \{0, 1\}^{|\setA_1|}\). Similarly, we denote the truncated vector for the toll-free arcs \(\xB \in \{0, 1\}^{|\setA_2|}\).
Once we fix \(\xA\), then the remaining program becomes a linear program,
which can be solved in polynomial time \citep{karmarkar1984}. As mentioned before, the single-commodity NPP is NP-hard \citep{roch2005},
therefore, the number of values of \(\xA\) that we need to enumerate cannot be polynomial (assuming P \(\neq\) NP).

Note that after \(\xA\) is fixed, we can decompose \cref{prog:1comm-1level} into two steps: (i)  Filling \(x\) with toll-free arcs \(\xB\) to form the shortest path; (ii) Finding the toll prices \(t\) that make \(x\) optimal.
If either step is infeasible, then \(\xA\) is infeasible in \cref{prog:1comm-1level}.
In step (i), the prices \(t\) are not required because the costs of the toll-free arcs do not depend on \(t\). The problem of
finding \(\xB\) becomes a minimum-cost flow problem with multiple sources and sinks, which can be solved by the following linear program:
\begin{equation*}
	\min_{\xB} \lbrace (\cx)_{\setA_2} \mid (\Ax)_{\setA_2} = b - (\Ax)_{\setA_1}, \xB \geq 0\rbrace.
\end{equation*}
Thus, once \(\xA\) is chosen, we know \(\xB\) immediately, and a path \(x\) emerges without the need of knowing \(t\).
It does not matter if there are multiple \(\xB\) for a single \(\xA\), since \(\xB\) does not affect the leader's revenue.
Having \(x\), the toll prices \(t\) can be computed in step (ii) with another linear program:
\begin{equation*}
	\max_{t,y}\lbrace  \xt  \mid  t \in \setT, \Ay - t \leq c, \by - \xt = \cx\rbrace.
\end{equation*}
From this perspective, the single-commodity NPP is a combinatorial problem,
which is less about the choice of the continuous prices, but more about the choice of the discrete paths (more specifically, the choice of \(\xA\)).
This observation will become more apparent once we introduce the conjugate model in~\cref{ssec:conjugate}.

\begin{example}\label{ex:single-comm}
	Consider a single-commodity NPP with the graph shown in \cref{fig:single-comm-network}.
	Dashed arcs are tolled arcs while solid arcs are toll-free. The numbers in the middle of the arcs
	represent their costs (arcs without numbers have no costs). In this graph, there are only two tolled arcs,
	hence we can compute the optimal solution by enumerating all 4 combinations of \(\xA = (x_1, x_2)\)
	(\(x_a\) corresponds to \(t_a, a \in \{1, 2\}\)):
	\begin{itemize}
		\item For \(\xA = (0, 0)\), the shortest path is \((o-d)\) with cost \(5\).
		      The leader receives no revenue.
		\item For \(\xA = (1, 0)\), the shortest path is \((o-u-v-p-d)\) with cost \(4 + t_1\).
		      Solving for \(t\) in step (ii) gives us \(t_1 = 1, t_2 \to \infty\). Thus, the leader's revenue~is~1.
		\item For \(\xA = (0, 1)\), the shortest path is \((o-v-p-q-d)\) with cost \(6 + t_2\).
		      However, the linear program in step (ii) is infeasible.
		\item For \(\xA = (1, 1)\), the shortest path is \((o-u-v-p-q-d)\) with cost \(2 + t_1 + t_2\).
		      The revenue is 3 with \(t_1 = 3, t_2 = 0\) (there are multiple solutions for \(t\),
		      all of them produce the same revenue).
	\end{itemize}
	Choosing the highest revenue from all 4 cases, we obtain the optimal value~of~3, corresponding to \(\xA = (1, 1)\) and the path
	\((o-u-v-p-q-d)\). Two remarks are in order from this example. First, in the case with \(\xA = (1, 0)\), the tolled arc 2 is not selected (\(x_2 = 0\)). Thus, this arc can be removed from the graph. We effectively achieved this by setting \(t_2\) to a very large number (basically to infinity). Second, in the case with \(\xA = (0, 1)\), the reason why step (ii) is infeasible is that the cost of this case (\(6 + t_2\)) is already larger than the cost of the first case (which is 5). Hence, no value \(t_2 \geq 0\) can make this path optimal. We call such a path \emph{bilevel infeasible}.  Bilevel feasibility is the topic of discussion in \cref{ssec:bf}.

	\begin{figure}
		\centering
		\scalebox{0.8}{\begin{tikzpicture}[scale=1.5]
    \tikzstyle{vertex}=[circle,draw,minimum size=20pt,inner sep=0pt]
    \tikzstyle{tolled}=[->, >=latex, dashed]
    \tikzstyle{tollfree}=[->, >=latex]

    \node[vertex] (o) at (0, 1) {$o$};
    \node[vertex] (u) at (0, 0) {$u$};
    \node[vertex] (v) at (1, 0) {$v$};
    \node[vertex] (p) at (2, 0) {$p$};
    \node[vertex] (q) at (3, 0) {$q$};
    \node[vertex] (d) at (3, 1) {$d$};

    \draw[tolled] (u) edge node[below]{$t_1$} (v) (p) edge node[below]{$t_2$} (q);
    \draw[tollfree] (o) edge (u) (v) edge node[above]{2} (p) (q) edge (d);
    \draw[tollfree] (o) edge node[above] {5} (d);
    \draw[tollfree] (o) edge node[above right]{4} (v);
    \draw[tollfree] (p) edge node[above left]{2} (d);

\end{tikzpicture}}
		\caption{Graph for \cref{ex:single-comm}.}
		\label{fig:single-comm-network}
	\end{figure}
\end{example}

\subsection{Single-Tolled-Arc Case}\label{ssec:single-toll}
We switch the focus on the other extreme case of the NPP, where there is a single tolled arc.
To stay away from triviality, we assume that there are several followers \(k \in \setK\),
each has its own origin \(\ok\) and destination \(\dk\) and optimizes its own shortest path \(\xk\) across the network.
The leader aims to maximize the sum of revenues from all followers.
The decision process of the leader can be thought of as the balancing act between
two strategies: sell high to a few, and sell low to everyone. We consider that every follower has the same demand \(\eta^k = 1\).
\footnote{All results presented in this paper can be extended to the non-unit demand case in a straightforward manner.}

Following the process described in \cref{ssec:single-comm}, the single-level reformulation of the single-tolled-arc NPP is:
\begin{subequations}
	\label[program]{prog:1toll-1level}
	\begin{align}
		\max_{t,x^k,y^k}\  & \sumk \txk       &  & \label[progeq]{eq:1toll-1level-obj}              \\
		\st                & t \in \setT,                                                           \\
		                   & \Axk = \bk,      &  & k \in \setK,                                     \\
		                   & \Ayk \leq c + t, &  & k \in \setK,                                     \\
		                   & \ctxk = \byk,    &  & k \in \setK,                                     \\
		                   & \xk \geq 0,      &  & k \in \setK. \label[progeq]{eq:1toll-1level-end}
	\end{align}
\end{subequations}
If we linearize this program as done in \cref{ssec:single-comm},
then the number of binary variables will be equal to the number of commodities \(|\setK|\), with one binary variable per commodity corresponding to
the tolled arc in \(\xk\). We set the index of this tolled arc to 0,
and let \(x_0 = (x^1_0, x^2_0, \ldots, x^{|\setK|}_0) \in \{0, 1\}^{|\setK|}\) be
the aggregated vector of the selections of the tolled arc from all commodities.
Since the number of binary variables is proportional to \(|\setK|\), is the single-tolled-arc NPP an NP-hard problem as well?

It turns out that the single-tolled-arc NPP can be solved in polynomial time. Understanding the algorithm used to solve this variant
is key to comprehend the asymmetry between \(|\setA_1|\) and \(|\setK|\). It is best to demonstrate this algorithm with an example.

\begin{example}\label{ex:single-toll}
	Consider the single-tolled-arc NPP described in \cref{fig:single-toll-network}.
	In this example, we only have one tolled arc which is shared by 3 different followers,
	each travels from \(\ok\) to \(\dk\) for \(k \in \setK = \{1, 2, 3\}\). Each follower can choose to use its own
	toll-free path with different cost (10, 4, and 3) or the tolled arc with cost \(t\).
	In fact, all instances of the single-tolled-arc NPP can be represented in this form. In \cref{ssec:single-comm},
	we already discussed that if \(\xkA\) is determined, so is \(\xkB\), and since \(|\setA_1| = 1\), for each follower,
	there are only 2 values of \(\xkA\) to consider: \(\xk_0 = 1\) and \(\xk_0 = 0\).

	Therefore, follower \(k\) will only use the tolled arc when \(t \leq c_k\) and use the toll-free path when \(t > c_k\).
	Now, we enumerate all 8 combinations of \(x_0 = (x^1_0, x^2_0, x^3_0) \in \{0, 1\}^3\).
	A quick observation tells us that we can never have \(x^1_0 = 0\) and \(x^2_0 = 1\) at the same time, since this implies
	\(t > 10\) and \(t \leq 4\) which is a contradiction. Generally, if \(c_{k_1} \geq c_{k_2}\) for some \(k_1, k_2 \in \setK\),
	then the combination \(x^{k_1}_0 = 0\) and \(x^{k_2}_0 = 1\) is infeasible. This property actually eliminates most of the combinations
	except for exactly \(|\setK| + 1 = 4\) combinations. Each combination corresponds to an interval of the price \(t\): \([0, 3]\), \((3, 4]\),
	\((4, 10]\), and \((10, \infty)\). These intervals are illustrated in the revenue plot in \cref{fig:single-toll-revplot}.
	The label of each interval denotes the set of commodities that use the tolled arc (\eg if \(t \in (3, 4]\),
	then the tolled arc is used by followers 1 and 2).

	From the revenue plot, we can conclude that the optimal revenue is 10, the optimal solution is \(t = 10\) and \((x^1_0, x^2_0, x^3_0) = (1, 0, 0)\).
	A remark regarding the 4 combinations is the following: each combination is characterized solely by the number of commodities that use the tolled arc, which
	happens to be 3, 2, 1, and 0 in this example (this explains the number \(|\setK| + 1\) above). This is the key difference between the single-commodity case
	and the single-tolled-arc case. In the single-commodity NPP, all (or most) combinations of \(\xA\) matter, whose number is \(2^{|\setA_1|}\) in total.
	In the single-tolled-arc NPP, only the number of commodities that use the tolled arc matters, so only \(|\setK| + 1\) cases need to be considered.

	\begin{figure}
		\centering
		\begin{subfigure}[b]{0.45\textwidth}
			\centering
			\scalebox{0.8}{\begin{tikzpicture}[xscale=1.5]
    \tikzstyle{vertex}=[circle,draw,minimum size=20pt,inner sep=0pt]
    \tikzstyle{tolled}=[->, >=latex, dashed]
    \tikzstyle{tollfree}=[->, >=latex]

    \node[vertex] (o1) at (0.6, 2) {$o^1$};
    \node[vertex] (d1) at (2.4, 2) {$d^1$};
    \node[vertex] (o2) at (0.3, 1) {$o^2$};
    \node[vertex] (d2) at (2.7, 1) {$d^2$};
    \node[vertex] (o3) at (0, 0) {$o^3$};
    \node[vertex] (d3) at (3, 0) {$d^3$};
    \node[vertex] (u) at (1, -1) {$u$};
    \node[vertex] (v) at (2, -1) {$v$};

    \node[inner sep=5pt,draw=none] (m1) at (3.5, 2) {};
    \node[inner sep=5pt,draw=none] (m2) at (3.5, 1) {};
    \node[inner sep=5pt,draw=none] (m3) at (3.5, 0) {};

    \draw[tollfree] (o1) edge node[above]{$c_1 = 10$} (d1);
    \draw[tollfree] (o2) edge node[above]{$c_2 = 4$} (d2);
    \draw[tollfree] (o3) edge node[above]{$c_3 = 3$} (d3);
    \draw[tollfree] (o1) edge (u) (o2) edge (u) (o3) edge (u);
    \draw[tollfree] (v) edge (d1) (v) edge (d2) (v) edge (d3);
    \draw[tolled] (u) edge node[above]{$t$} (v);
\end{tikzpicture}}
			\subcaption{Graph}
			\label{fig:single-toll-network}
		\end{subfigure}
		\hfill
		\begin{subfigure}[b]{0.52\textwidth}
			\centering
			\scalebox{0.8}{\begin{tikzpicture}
    \tikzstyle{label}=[inner sep=5pt]
    \tikzstyle{soliddot}=[only marks,mark=*]
    \tikzstyle{hollowdot}=[fill=white,only marks,mark=*]
    \begin{axis}[
            scale = 0.65,
            axis lines = left,
            ylabel = Leader's revenue,
            xlabel = Toll price $t$,
            xtick distance = 2,
        ]
        \addplot[domain=0:3] {3*x} node[label, pos=0.2, right] {$\{1, 2, 3\}$};
        \addplot[domain=3:4] {2*x} node[label, pos=0.5, right] {$\{1, 2\}$};
        \addplot[domain=4:10] {x} node[label, pos=0.5, below right] {$\{1\}$};
        \addplot[domain=10:12] {0} node[label, pos=0.5, above] {$\varnothing$};
        \draw[dotted] (axis cs:3,9) -- (axis cs:3,6);
        \draw[dotted] (axis cs:4,8) -- (axis cs:4,4);
        \draw[dotted] (axis cs:10,10) -- (axis cs:10,0);
        \addplot[hollowdot] coordinates{(3,6)(4,4)(10,0)};
        \addplot[soliddot] coordinates{(3,9)(4,8)(10,10)};
    \end{axis}
\end{tikzpicture}}
			\subcaption{Revenue plot}
			\label{fig:single-toll-revplot}
		\end{subfigure}
		\caption{Illustrations for \cref{ex:single-toll}.}
	\end{figure}
\end{example}

After having our intuition established, we proceed to formalize this process in \cref{alg:single-toll}.
The first part of the algorithm (lines 1-5) calculates the maximum price that follower \(k\) will remain using the tolled arc.
If \(t \leq c_k\), then the follower will use the tolled arc and vice-versa. In some cases, \(c_k = 0\) due to \(\overline{c}_k = \underline{c}_k\),
then the tolled path is never as good as the toll-free path and no values of \(t \geq 0\) will persuade this follower to use the tolled arc
(another example of \emph{bilevel infeasibility}). Next, we sort all \(c_k\) (line 6) so that when \(t = c_{k_1}\), then only \(k_1\) uses the tolled arc,
when \(t = c_{k_2}\), then \(k_1\) and \(k_2\) use the tolled arc, and so on. Then, we enumerate all cases (lines 7-11) according to the number of commodities
that use the tolled arc, denoted as \(w\) (line 8). The revenue \(r_w\) of the case \(w\) is the product of the number of commodities \(w\)
and the maximum price \(c_{k_w}\) such that \(w\) commodities still use the tolled arc (line 9). Finally, we collect the maximum of all \(r_w\) (line 10)
and return it (line 12). Note that in \cref{alg:single-toll}, we skip the case \(w = 0\) since it does not produce any revenue, but theoretically,
this case still exists and it will become relevant in the general case.

\begin{algorithm}
	\caption{Single-tolled-arc NPP in polynomial time \citep{labbe2021}}
	\label{alg:single-toll}

	\begin{algorithmic}[1] \footnotesize
		\Require{The graph \(G = (\setV, \setA)\) with a single tolled arc,
			set of commodities~\(\setK\) with their O-D pairs \(o^k, d^k\).}
		\Ensure{The optimal leader's revenue.}
		\ForAll{\(k \in \setK\)}
		\State \(\overline{c}_k \gets\) cost of the shortest path when \(t \to \infty\)
		\State \(\underline{c}_k \gets\) cost of the shortest path when \(t = 0\)
		\State \(c_k \gets \overline{c}_k - \underline{c}_k\)
		\EndFor
		\State Sort \(c_k\) in descending order \(c_{k_1}, c_{k_2}, \ldots, c_{k_{|\setK|}}\)
		\State \(R \gets 0\) \Comment{Maximum revenue}
		\For{\(w \gets 1, 2, \ldots, |\setK|\)}
		\State \(r_w = w c_{k_w}\)
		\State \(R \gets \max\{R, r_w\}\)
		\EndFor
		\State \Return \(R\)
	\end{algorithmic}
\end{algorithm}

\subsection{Reaction Plot}\label{ssec:reaction}
In the remaining of \cref{sec:asymmetry}, we will consider the general case where there are multiple commodities and multiple tolled arcs.
Because the single-commodity case can be reduced to the general case, the general case is also an NP-hard problem.
Extending \cref{prog:1comm-1level} to the general case results in a number of binary variables equal to \(|\setA_1||K|\),
consisting of \(\xk_a\) for each \(a \in \setA_1\) and \(k \in \setK\). However, as shown in \cref{ssec:single-toll},
this number is not a good indicator for the complexity.

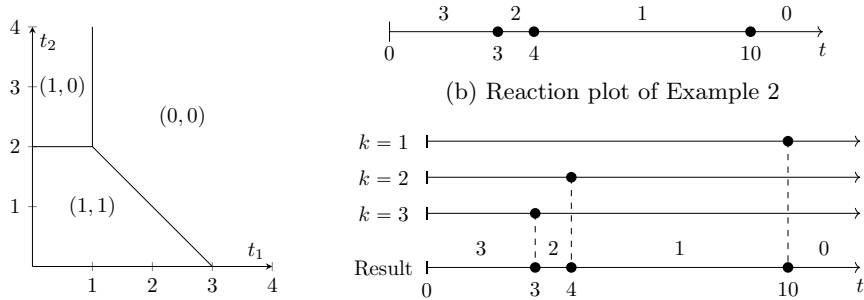
\begin{figure}
	\begin{minipage}[b][][b]{0.38\textwidth}
		\begin{subfigure}[b]{\textwidth}
			\centering
			\scalebox{0.8}{\begin{tikzpicture}
    \begin{axis}[
            scale=0.7,
            axis equal image,
            axis lines = center,
            xlabel = $t_1$,
            ylabel = $t_2$,
            xmin=0,
            xmax=4,
            ymin=0,
            ymax=4,
        ]

        \draw (3,0) -- (1,2) -- (1,4);
        \draw (0,2) -- (1,2);
        \node at (1,1) {$(1,1)$};
        \node at (0.5, 3) {$(1,0)$};
        \node at (2.5, 2.5) {$(0,0)$};
    \end{axis}
\end{tikzpicture}}
			\subcaption{Reaction plot of \cref{ex:single-comm}}
			\label{fig:single-comm-reaction}
		\end{subfigure}
	\end{minipage}
	\hfill
	\begin{minipage}[b][][b]{0.6\textwidth}
		\centering
		\begin{subfigure}[b]{\textwidth}
			\centering
			\scalebox{0.8}{\begin{tikzpicture}[scale=0.6]
    \tikzstyle{soliddot}=[circle,fill=black,inner sep=0pt,minimum size=5pt]
    \tikzstyle{belowlabel}=[below=4pt of #1,inner sep=0pt]
    \tikzstyle{abovelabel}=[above=2pt of #1,inner sep=0pt]
    \draw[->] (0,0) -- (12,0) node[below=2pt] {$t$};
    \draw (0,8pt) -- (0,-8pt);
    \node[soliddot,fill=none] at (0,0) (t0) {} node[belowlabel=t0] {0};
    \node[soliddot] (t3) at (3, 0) {} node[belowlabel=t3] {3};
    \node[soliddot] (t4) at (4, 0) {} node[belowlabel=t4] {4};
    \node[soliddot] (t10) at (10, 0) {} node[belowlabel=t10] {10};

    \node (w3) at (1.5, 0) {} node[abovelabel=w3] {3};
    \node (w2) at (3.5, 0) {} node[abovelabel=w2] {2};
    \node (w1) at (7, 0) {} node[abovelabel=w1] {1};
    \node (w0) at (11, 0) {} node[abovelabel=w0] {0};
\end{tikzpicture}}
			\subcaption{Reaction plot of \cref{ex:single-toll}}
			\label{fig:single-toll-reaction}
		\end{subfigure}
		\par\medskip
		\begin{subfigure}[b]{\textwidth}
			\centering
			\scalebox{0.8}{\begin{tikzpicture}[scale=0.6]
    \tikzstyle{soliddot}=[circle,fill=black,inner sep=0pt,minimum size=5pt]
    \tikzstyle{origin}=[fill=black,inner sep=0pt,minimum height=6pt, minimum width=0.6pt]
    \tikzstyle{belowlabel}=[below=4pt of #1,inner sep=0pt]
    \tikzstyle{abovelabel}=[above=2pt of #1,inner sep=0pt]

    \def\vsep{1}

    \draw[->] (0,0) node[origin] (t0) {} -- (12,0) node[below=2pt] {$t$};
    \node[belowlabel=t0] {0};
    \node[soliddot] (t3) at (3, 0) {} node[belowlabel=t3] {3};
    \node[soliddot] (t4) at (4, 0) {} node[belowlabel=t4] {4};
    \node[soliddot] (t10) at (10, 0) {} node[belowlabel=t10] {10};

    \node (w3) at (1.5, 0) {} node[abovelabel=w3] {3};
    \node (w2) at (3.5, 0) {} node[abovelabel=w2] {2};
    \node (w1) at (7, 0) {} node[abovelabel=w1] {1};
    \node (w0) at (11, 0) {} node[abovelabel=w0] {0};

    \draw[->] (0,0.5+\vsep) node[origin] {} -- ++ (12,0);
    \draw[->] (0,0.5+2*\vsep) node[origin] {} -- ++ (12,0);
    \draw[->] (0,0.5+3*\vsep) node[origin] {} -- ++ (12,0);

    \node[soliddot] at (3, 0.5+\vsep) {};
    \node[soliddot] at (4, 0.5+2*\vsep) {};
    \node[soliddot] at (10, 0.5+3*\vsep) {};

    \draw[dashed] (3,0.5+\vsep) -- (3,0);
    \draw[dashed] (4,0.5+2*\vsep) -- (4,0);
    \draw[dashed] (10,0.5+3*\vsep) -- (10,0);

    \node at (-1.2, 0) {Result};
    \node at (-1.2, 0.5+\vsep) {\(k = 3\)};
    \node at (-1.2, 0.5+2*\vsep) {\(k = 2\)};
    \node at (-1.2, 0.5+3*\vsep) {\(k = 1\)};
\end{tikzpicture}}
			\subcaption{Composition of the reaction plot in \cref{ex:single-toll}.}
			\label{fig:single-toll-composition}
		\end{subfigure}
	\end{minipage}
	\caption{Examples of reaction plots.}
\end{figure}

To extend the result presented in \cref{ssec:single-toll}, we employ an illustration tool called \emph{reaction plot}.
In the reaction plot, the price of each tolled arc is represented by an axis (\(n\) tolled arcs means \(n\) axes).
Then, for each value of \(t\), we plot and group the reaction of the followers given \(t\), \ie the shortest path in the NPP.
The result is a partition of the values of \(t\) into separate sections corresponding to different reactions of the followers.
The reaction plot of \cref{ex:single-toll} is shown in \cref{fig:single-toll-reaction}. The axis represents the price \(t\) of the tolled arc.
The numbers below the axis are the meeting points of the intervals (which are \(c_k\)), and the numbers above the axis are the number of commodities that use
the tolled arc for each interval (which is \(w\) in \cref{alg:single-toll}). The reaction plot of \cref{ex:single-comm}, displayed in \cref{fig:single-comm-reaction},
is a 2-dimensional plot instead, since there are 2 tolled arcs in this case. The label of each region is \(\xA\).
Note that there is no region for \(\xA = (0, 1)\) since, as we remarked in \cref{ssec:single-comm},
there is no value of \(t \geq 0\) that makes \(\xA = (0, 1)\) an optimal follower's solution (bilevel infeasible).

We have explored the reaction plots of the single-commodity and the single-tolled-arc case. Now, the question is:
How do we draw the reaction plot of the general case? One could solve the followers' problems for each value of \(t\), then
group the reactions as we did in \cref{fig:single-comm-reaction,fig:single-toll-reaction}. Instead, we will
compose the reaction plot of the general case from the plots of multiple single-commodity subproblems, one subproblem per commodity.
Consider the reaction plot in~\cref{fig:single-toll-reaction}. This plot is actually a composition of three
individual reaction plots, superimposed upon each other as described in~\cref{fig:single-toll-composition}.
The reaction plot of a single commodity \(k \in \setK\) in this case contains only 2 intervals: \([0, c_k]\) and \((c_k, \infty)\)
(corresponding to \(\xk_0 = 1\) and \(\xk_0 = 0\), respectively). We are allowed to stack up the reaction plots because
once \(t\) is set, each follower solves its own follower's problem separately, thus their reactions are independent.
This principle still applies if we have multiple tolled arcs, as shown in the next example.

\begin{example}\label{ex:multi-comm}
	Consider a multi-commodity NPP with graph in \cref{fig:multi-comm-network}. There are two commodities
	travelling from \(o^k\) to \(d^k\), \(k \in \{1, 2\}\). There are also two tolled arcs with prices \(t_1\) and \(t_2\),
	shared by both commodities. With respect to each follower, only half the graph is relevant.
	\Cref{fig:multi-comm-series-network,fig:multi-comm-parallel-network} are the subgraphs in the perspectives
	of followers 1 and 2, respectively. Irrelevant nodes are drawn with dashes border while irrelevant arcs are hidden.
	The reaction plots of individual commodities are displayed in \cref{fig:multi-comm-series-reaction,fig:multi-comm-parallel-reaction}.
	By stacking \cref{fig:multi-comm-series-reaction} on top of \cref{fig:multi-comm-parallel-reaction},
	we have the overall reaction plot in \cref{fig:multi-comm-reaction}. Although each individual reaction plot has 4 regions,
	the composed plot only has \(8\) regions in total rather than \(4 \times 4 = 16\).

	\begin{figure}
		\centering

		\begin{subfigure}[b]{0.32\textwidth}
			\centering
			\scalebox{0.8}{\begin{tikzpicture}[scale=1.4]
    \tikzstyle{vertex}=[circle,draw,minimum size=20pt,inner sep=0pt]
    \tikzstyle{tolled}=[->, >=latex, dashed]
    \tikzstyle{tollfree}=[->, >=latex]

    \node[vertex] (o1) at (1, 1) {$o^1$};
    \node[vertex] (o2) at (1, -1) {$o^2$};
    \node[vertex] (u) at (1, 0) {$u$};
    \node[vertex] (v) at (2, 0) {$v$};
    \node[vertex] (p) at (3, 0) {$p$};
    \node[vertex] (q) at (4, 0) {$q$};
    \node[vertex] (d1) at (4, 1) {$d^1$};
    \node[vertex] (d2) at (4, -1) {$d^2$};

    \draw[tollfree] (v) edge node[above]{3} (p);
    \draw[tolled] (u) edge node[below]{$t_1$} (v) (p) edge node[below]{$t_2$} (q);
    
    \draw[tollfree] (o1) edge (u) (q) edge (d1);
    \draw[tollfree] (o1) edge node[above] {12} (d1);
    \draw[tollfree] (o1) edge node[above right]{5} (v);
    \draw[tollfree] (p) edge node[above left]{6} (d1);

    \draw[tollfree] (o2) edge (u) (q) edge (d2);
    \draw[tollfree] (o2) edge node[below] {11} (d2);
    \draw[tollfree] (o2) edge node[below]{5} (p);
    \draw[tollfree] (v) edge node[below]{6} (d2);

\end{tikzpicture}}
			\subcaption{Overall graph}
			\label{fig:multi-comm-network}
		\end{subfigure}
		\hfill
		\begin{subfigure}[b]{0.32\textwidth}
			\centering
			\scalebox{0.8}{\begin{tikzpicture}[scale=1.4]
    \tikzstyle{vertex}=[circle,draw,minimum size=20pt,inner sep=0pt]
    \tikzstyle{tolled}=[->, >=latex, dashed]
    \tikzstyle{tollfree}=[->, >=latex]

    \node[vertex] (o1) at (1, 1) {$o^1$};
    \node[vertex, dashed] (o2) at (1, -1) {$o^2$};
    \node[vertex] (u) at (1, 0) {$u$};
    \node[vertex] (v) at (2, 0) {$v$};
    \node[vertex] (p) at (3, 0) {$p$};
    \node[vertex] (q) at (4, 0) {$q$};
    \node[vertex] (d1) at (4, 1) {$d^1$};
    \node[vertex, dashed] (d2) at (4, -1) {$d^2$};

    \draw[tollfree] (v) edge node[above]{3} (p);
    \draw[tolled] (u) edge node[below]{$t_1$} (v) (p) edge node[below]{$t_2$} (q);
    
    \draw[tollfree] (o1) edge (u) (q) edge (d1);
    \draw[tollfree] (o1) edge node[above] {12} (d1);
    \draw[tollfree] (o1) edge node[above right]{5} (v);
    \draw[tollfree] (p) edge node[above left]{6} (d1);

    \draw[opacity=0] (o2) edge node[below] {\phantom{11}} (d2);

\end{tikzpicture}}
			\subcaption{Subgraph for \(k = 1\)}
			\label{fig:multi-comm-series-network}
		\end{subfigure}
		\hfill
		\begin{subfigure}[b]{0.32\textwidth}
			\centering
			\scalebox{0.8}{\begin{tikzpicture}[scale=1.4]
    \tikzstyle{vertex}=[circle,draw,minimum size=20pt,inner sep=0pt]
    \tikzstyle{tolled}=[->, >=latex, dashed]
    \tikzstyle{tollfree}=[->, >=latex]

    \node[vertex, dashed] (o1) at (1, 1) {$o^1$};
    \node[vertex] (o2) at (1, -1) {$o^2$};
    \node[vertex] (u) at (1, 0) {$u$};
    \node[vertex] (v) at (2, 0) {$v$};
    \node[vertex] (p) at (3, 0) {$p$};
    \node[vertex] (q) at (4, 0) {$q$};
    \node[vertex, dashed] (d1) at (4, 1) {$d^1$};
    \node[vertex] (d2) at (4, -1) {$d^2$};

    \draw[tollfree] (v) edge node[above]{3} (p);
    \draw[tolled] (u) edge node[below]{$t_1$} (v) (p) edge node[below]{$t_2$} (q);
    

    \draw[tollfree] (o2) edge (u) (q) edge (d2);
    \draw[tollfree] (o2) edge node[below] {11} (d2);
    \draw[tollfree] (o2) edge node[below]{5} (p);
    \draw[tollfree] (v) edge node[below]{6} (d2);

\end{tikzpicture}}
			\subcaption{Subgraph for \(k = 2\)}
			\label{fig:multi-comm-parallel-network}
		\end{subfigure}
		\par\medskip

		\begin{subfigure}[b]{0.32\textwidth}
			\centering
			\scalebox{0.8}{\begin{tikzpicture}
    \begin{axis}[
            scale=0.7,
            axis equal image,
            axis lines = center,
            xlabel = $t_1$,
            ylabel = $t_2$,
            xmin=0,
            xmax=8,
            ymin=0,
            ymax=8,
        ]
        \draw (5,0) -- (5,4) -- (8,4);
        \draw (3,8) -- (3,6) -- (0,6);
        \draw (3,6) -- (5,4);

        \draw (2,0) -- (2,3) -- (0,3);
        \draw (5,8) -- (5,6) -- (8,6);
        \draw (2,3) -- (5,6);

        \node at (1,1.5) {$(2,2)$};
        \node at (1.7,4.7) {$(2,1)$};
        \node at (2.0,7) {$(2,0)$};
        \node at (3.5,2) {$(1,2)$};
        \node at (4,7) {$(1,0)$};
        \node at (6.5,2) {$(0,2)$};
        \node at (6.5,5) {$(0,1)$};
        \node at (6.5,7) {$(0,0)$};
    \end{axis}
\end{tikzpicture}}
			\subcaption{Composed reaction plot}
			\label{fig:multi-comm-reaction}
		\end{subfigure}
		\hfill
		\begin{subfigure}[b]{0.32\textwidth}
			\centering
			\scalebox{0.8}{\begin{tikzpicture}
    \begin{axis}[
        scale=0.7,
            axis equal image,
            axis lines = center,
            xlabel = $t_1$,
            ylabel = $t_2$,
            xmin=0,
            xmax=8,
            ymin=0,
            ymax=8,
        ]
        \draw (5,0) -- (5,4) -- (8,4);
        \draw (3,8) -- (3,6) -- (0,6);
        \draw (3,6) -- (5,4);
        \node at (2.5,3) {$(1,1)$};
        \node at (2.0,7) {$(1,0)$};
        \node at (6.5,2) {$(0,1)$};
        \node at (5.5,6) {$(0,0)$};
    \end{axis}
\end{tikzpicture}}
			\subcaption{Reaction plot for \(k = 1\)}
			\label{fig:multi-comm-series-reaction}
		\end{subfigure}
		\hfill
		\begin{subfigure}[b]{0.32\textwidth}
			\centering
			\scalebox{0.8}{\begin{tikzpicture}
    \begin{axis}[
        scale=0.7,
            axis equal image,
            axis lines = center,
            xlabel = $t_1$,
            ylabel = $t_2$,
            xmin=0,
            xmax=8,
            ymin=0,
            ymax=8,
        ]
        \draw (2,0) -- (2,3) -- (0,3);
        \draw (5,8) -- (5,6) -- (8,6);
        \draw (2,3) -- (5,6);
        \node at (1,1.5) {$(1,1)$};
        \node at (2,6) {$(1,0)$};
        \node at (5.5,2.5) {$(0,1)$};
        \node at (6.5,7) {$(0,0)$};
    \end{axis}
\end{tikzpicture}}
			\subcaption{Reaction plot for \(k = 2\)}
			\label{fig:multi-comm-parallel-reaction}
		\end{subfigure}

		\caption{Composition of the reaction plot in \cref{ex:multi-comm}.}
	\end{figure}
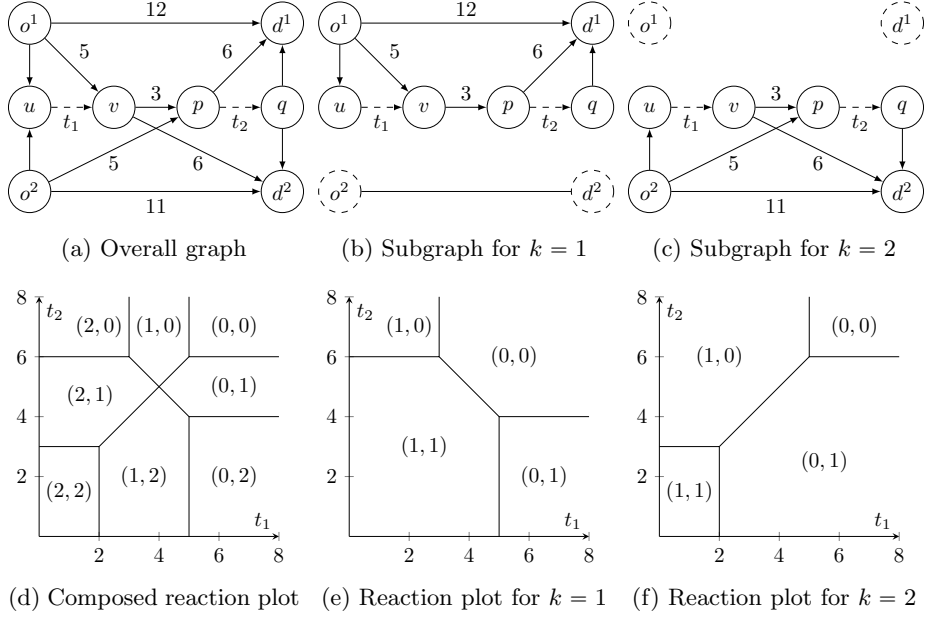

	The labels of the regions in \cref{fig:multi-comm-reaction} are in the format \((w_1, w_2)\), where
	\(w_1 = x_1^1 + x_1^2\) is the number of commodities that uses the first tolled arc,
	while \(w_2 = x_2^1 + x_2^2\) is the same but for the second tolled arc. Similar to the single-tolled-arc case in
	\cref{ssec:single-toll},	when stacking the reaction plots, some combinations are prohibited due to
	the geometry of the individual plots. For example, we only have \(\xA^1 = (1, 1)\) in combination with \(\xA^2 = (0, 1)\),
	but not \(\xA^1 = (0, 1)\) and \(\xA^2 = (1, 1)\). Thus, the label \((1, 2)\) in \cref{fig:multi-comm-reaction} always
	indicates the former combination, not the latter. Using this labeling system, intuitively, we can only have a maximum of
	\((|\setK| + 1)^{|\setA_1|}\) different regions in the reaction plot.
	The maximum number of regions in this particular example is \(3^2 = 9\), but we can have less,
	\eg in \cref{fig:multi-comm-reaction}, we only have 8 regions (the label \((1, 1)\) is missing,
	to which we will come back in \cref{ex:multi-comm-calc}).

	A side note is that the graph in this example is specifically designed as a showcase of two different structures in the case
	of two tolled arcs. The graph in \cref{fig:multi-comm-series-network} has the \emph{series} structure,
	whose reaction plot (\cref{fig:multi-comm-series-reaction}) contains a diagonal line from the top-left corner to the bottom-right.
	On the other hand, \cref{fig:multi-comm-parallel-network} represents the \emph{parallel} structure (remove the arc \(v-p\)
	for more clarity), whose reaction plot (\cref{fig:multi-comm-parallel-reaction}) also has a diagonal line
	but in the other direction (bottom-left to top-right). The overall graph demonstrates the complex nature of the NPP,
	where a pair of tolled arcs can be in series for one commodity, while in parallel for another.
\end{example}

\subsection{Conjugate Model}\label{ssec:conjugate}
The labeling system \((w_1, w_2, \ldots, w_{|\setA_1|})\) from \cref{ex:multi-comm} gives us a hint:
the number of discrete cases that we need to consider is only \((|\setK| + 1)^{|\setA_1|}\), rather than \(2^{|\setA_1||K|}\)
different combinations of the binary variables \(\xk_a\). The next step is to formalize our intuition and prove \cref{thm:multi-comm-poly}.
For this task, we developed a new reformulation of the NPP called the \emph{conjugate model}.

We start with the bilevel program of the multi-commodity case:
\begin{equation}
	\label[program]{prog:npp}
	\max_{t,\xk} \Bigl\{ \sumk \txk \mid t \in \setT, \xk \in \funcR^k(t) \ \ \forall k \in \setK \Bigr\},
\end{equation}
where \(\funcR^k(t)\) is the reaction set of follower \(k \in \setK\):
\begin{equation*}
	\funcR^k(t) = \argmin_{\xk}\Bigl\{  \ctxk  \mid \Axk = \bk, \xk \geq 0 \Bigr\}.
\end{equation*}
Once \(t\) is set by the leader, all followers' problems are independent. Thus, we can combine all of these problems into
a single aggregated problem:
\begin{equation*}
	\overline{\funcR}(t) = \argmin_{\xk} \Bigl\{   \sumk \ctxk  \mid  \Axk = \bk \ \ \forall  k \in \setK, \xk \geq 0 \ \ \forall k \in \setK  \Bigr\}.
\end{equation*}
The aggregated reaction set \(\overline{\funcR}(t)\) is equivalent to the product of all the individual reaction sets: $\overline{\funcR}(t) = \prod_{k\in\setK} \funcR^k(t).$

Next, we introduce a new variable \(w = \sumk \xkA\). The vector \(w\) has the same meaning as \((w_1, w_2, \ldots, w_{|\setA_1|})\)
used in the labeling system, where \(w_a\) is the number of commodities that use the arc \(a \in \setA_1\).
However, we do not replace \(\sumk \xkA\) by \(w\) with an equality, but rather with an inequality:
\begin{subequations}
	\label[program]{prog:aggregated}
	\begin{align}
		\min_{w, \xk}\  & \sumk \left(\cxk\right) + \tw \label[progeq]{eq:aggregated-obj}                                                     \\
		\st             & \sumk \xkA \leq w, \label[progeq]{eq:aggregated-xw}                                                                 \\
		                & \Axk = \bk,                                                     &  & k \in \setK,                                   \\
		                & w \geq 0,                                                                                                           \\
		                & \xk \geq 0,                                                     &  & k \in \setK. \label[progeq]{eq:aggregated-end}
	\end{align}
\end{subequations}
As an abuse of notation, we used \(t = t_{\setA_1}\) in the objective function~\cref{eq:aggregated-obj},
which is the truncated version of the original~\(t\) (the version of \(t\) is implied by the context).
Note that the objective function minimizes \(\tw\), hence the inequality~\cref{eq:aggregated-xw} is always active whenever \(t > 0\).
The dual of \cref{prog:aggregated} is:
\begin{equation}
	\label[program]{prog:aggregated-dual}
	\max_{\yk} \Bigl\{  \sumk \byk \mid  \Ayk \leq c + t \ \ \forall k \in \setK \Bigr\}.
\end{equation}

Given \(w \geq 0\), we introduce the \emph{conjugate follower formulation}, defined as:
\begin{equation}
	\label[program]{prog:conjugate-follower}
	\max_{t, \yk}\Bigl\{    \sumk \left(\byk\right) - \wt  \mid \Ayk - t \leq c \ \ \forall k \in \setK,  t \geq 0\Bigr\}.
\end{equation}

\Cref{prog:conjugate-follower} is similar to \cref{prog:aggregated-dual}, except for the following changes: (i) \Cref{prog:conjugate-follower} is parameterized by \(w\), while \cref{prog:aggregated-dual} is parameterized by \(t\); (ii) The toll prices \(t\) become variables, and \(t \geq 0\) is added as a constraint; (iii) The term \(-\wt\) is added to the objective function of \cref{prog:conjugate-follower}.

The dual of \cref{prog:conjugate-follower} is:
\begin{equation}
	\label[program]{prog:conjugate-dual}
	\min_{\xk}\Bigl\{   \sumk \cxk \mid  \sumk \xkA \leq w,  \Axk = \bk \ \ \forall k \in \setK, \xk \geq 0 \ \  \forall k \in \setK\Bigr\}.
\end{equation}

Let \(\funcT(w)\) be the set of \(t\) that are optimal to \cref{prog:conjugate-follower} given \(w \geq 0\).
We define the \emph{conjugate bilevel formulation} as:
\begin{equation}
	\label[program]{prog:conjugate-leader}
	\max_{w,t}\lbrace  \wt  \mid w \geq 0,  t \in \funcT(w)\rbrace.
\end{equation}

\begin{proposition}\label{prop:conjugate-equivalence}
	\Cref{prog:npp,prog:conjugate-leader} are equivalent, in the sense that their optimal objective values are equal.
\end{proposition}
\begin{proof}
	We use the KKT conditions together with strong duality on \cref{prog:aggregated,prog:aggregated-dual}
	to convert \cref{prog:npp} to the following single-level reformulation:
	\begin{subequations}
		\label[program]{prog:aggregated-1level}
		\begin{align}
			\max_{t,w,x,y}\  & \tw                                                          \\
			\st              & \sumk \xkA \leq w,                                           \\
			                 & \Axk = \bk,                                &  & k \in \setK, \\
			                 & \Ayk \leq c + t,                           &  & k \in \setK, \\
			                 & \sumk \left(\cxk\right) + \tw = \sumk\byk,                   \\
			                 & t \geq 0,                                                    \\
			                 & w \geq 0,                                                    \\
			                 & \xk \geq 0,                                &  & k \in \setK.
		\end{align}
	\end{subequations}
	Applying the same technique to \cref{prog:conjugate-follower,prog:conjugate-dual} produces the same single-level reformulation
	for \cref{prog:conjugate-leader}. Therefore, the optimal objective values of \cref{prog:npp,prog:conjugate-leader} are the same.
\end{proof}

The \emph{conjugate model} provides a different perspective to NPP, where the leader controls \(w\) in contrast to \(t\) in the original model.
Looking at \cref{prog:conjugate-dual}, the role of \(w\) is that of the capacities of the tolled arcs.
\Cref{prog:conjugate-dual} is similar to the shortest path problem, but with capacities limited to \(w\).
Note that this program is not totally unimodular, thus there is no guarantee that \(\xk\) will be binary even if \(w\) is integral.

Hereafter, we refer to \(\funcT(w)\) as the \emph{action set} and reserve the term \emph{reaction set} for \(\funcR(t)\).
Similarly, we call \(t\) the \emph{action}, \(x\) the \emph{reaction}, and \(w\) the \emph{reduced reaction}.

\begin{proof}[Proof of \cref{thm:multi-comm-poly}]
	\Cref{prop:conjugate-equivalence} implies that we can use \cref{prog:conjugate-leader} to find the optimal revenue of the NPP.
	For each commodity \(k \in \setK\), we know that there is always an optimal reaction \(\xk\) that is binary (since the followers' problems are shortest path problems),
	and since \(w = \sumk \xkA\), it is sufficient to enumerate \(w\) in the set \(\{0, 1, \ldots, |\setK|\}^{|\setA_1|}\).
	For each \(w\), we solve \cref{prog:conjugate-follower} in polynomial time to obtain \(t\),
	which can be multiplied with \(w\) to compute the revenue. The total number of cases that we need to enumerate is \((|\setK| + 1)^{|\setA_1|}\).
	Thus, if the number of tolled arcs \(|\setA_1|\) is fixed, this process solves the multi-commodity NPP in polynomial time.
\end{proof}

Although we can solve the NPP in polynomial time given that \(|\setA_1|\) is fixed, the approach described in the proof of \cref{thm:multi-comm-poly} is not of practical interest.
For \(|\setA_1| = 2\), the complexity of the algorithm is in \(O(|\setK|^2)\), for \(|\setA_1| = 10\), it is in \(O(|\setK|^{10})\), and so on.
What \cref{thm:multi-comm-poly} emphasizes is that there exists an asymmetry in the complexity,
which we can exploit to tackle problems with a high number of commodities. This will be attempted in the next section, with the help of the conjugate model.

\begin{example}\label{ex:multi-comm-calc}
	In this example, we compute the optimal revenue of the NPP described in \cref{ex:multi-comm} using the algorithm in
	the proof of \cref{thm:multi-comm-poly}. Since \(|\setK| = 2\) and \(|\setA_1| = 2\), we have 9 possible values for \(w\),
	which are listed in \cref{tab:multi-comm-calc}.
	The optimal solution is highlighted in bold. The optimal revenue is 14, corresponding to \((w_1, w_2) = (1, 2)\) and \((t_1, t_2) = (4, 5)\).
	The case \(w = (1, 1)\) is special, since \cref{prog:conjugate-dual} returns a fractional \(\xk_a\) reaction instead of a binary result.
	Using \cref{fig:multi-comm-reaction} as a reference, this case is not represented by a full-dimensional region in the reaction plot, but by just a single
	point at \(t = (4, 5)\). We call this degenerate case \emph{weakly bilevel feasible}, as opposed to the other 8 cases which are \emph{strongly bilevel feasible}.
	Note that weakly bilevel feasible cases still appear in the reaction plot, but they do not manifest into full-dimensional regions. This contrasts with bilevel infeasible cases which
	are completely absent from the reaction plot. Weak and strong bilevel feasibility are discussed in~\cref{ssec:strongbf}.

	\begin{table}
		\caption{Enumeration of all cases in \cref{ex:multi-comm}.}
		\label{tab:multi-comm-calc}
		\centering \footnotesize
		\begin{tabular}{ccccccccc}
			\toprule
			        &         & \multicolumn{2}{c}{Actions} & \multicolumn{4}{c}{Reactions}                 & Revenue \\
			\cmidrule(lr){3-4}\cmidrule(lr){5-8}\cmidrule(lr){9-9}
			\(w_1\) & \(w_2\) & \(t_1\)    & \(t_2\)        & \(x^1_1\) & \(x^1_2\) & \(x^2_1\) & \(x^2_2\) & \(\wt\) \\
			\midrule
			0       & 0       & \(\infty\) & \(\infty\)     & 0         & 0         & 0         & 0         & 0       \\
			1       & 0       & 5          & \(\infty\)     & 0         & 0         & 1         & 0         & 5       \\
			2       & 0       & 3          & \(\infty\)     & 1         & 0         & 1         & 0         & 6       \\
			0       & 1       & \(\infty\) & 6              & 0         & 0         & 0         & 1         & 6       \\
			1       & 1       & 4          & 5              & 0.5       & 0.5       & 0.5       & 0.5       & 9       \\
			2       & 1       & 4          & 5              & 1         & 1         & 1         & 0         & 13      \\
			0       & 2       & \(\infty\) & 4              & 0         & 1         & 0         & 1         & 8       \\
			\bf 1   & \bf 2   & \bf 4      & \bf 5          & \bf 1     & \bf 1     & \bf 0     & \bf 1     & \bf 14  \\
			2       & 2       & 2          & 3              & 1         & 1         & 1         & 1         & 10      \\
			\bottomrule
		\end{tabular}
	\end{table}
\end{example}

\section{Strong Bilevel Feasibility}\label{sec:strongbf}
\Cref{sec:asymmetry} tells us that the number of reduced reactions \(w\) that we need to enumerate is \((|\setK| + 1)^{|\setA_1|}\).
This number is far smaller than the number of all compositions of integral values of \(\xkA\) which is \(2^{|\setA_1||\setK|}\).
Since \(w\) is defined by \(w = \sumk \xkA\), we can associate a value of \(w\) to a composition of \(\xkA\),
while a value of \(w\) can be associated to many compositions of \(\xkA\).
This raises several questions: How do we map \((|\setK| + 1)^{|\setA_1|}\) values of \(w\)  back into  \(2^{|\setA_1||\setK|}\) compositions of \(\xkA\)?
Can all compositions of \(\xkA\) be mapped from \(w\), or are there some compositions that we can eliminate?
What are the properties that a composition needs to satisfy to not be eliminated and how do we test them?

To answer these questions, we further explore the notion of the reduced reaction \(w\),
discuss the duality between \(t\) and \(w\), from which we derive the concept of \emph{strong bilevel feasibility}.
Strong bilevel feasibility is the property that will eliminate most of \(\xkA\) and bring the number \(2^{|\setA_1||\setK|}\) down to \((|\setK| + 1)^{|\setA_1|}\).
We return to the single-commodity case to define three concepts: \emph{bilevel feasibility} (\cref{ssec:bf}),
\emph{action-reaction duality} (\cref{ssec:duality}), and \emph{strong bilevel feasibility} (\cref{ssec:strongbf}).
Finally, we extend these concepts to the multi-commodity case via the composition of individual commodities (\cref{ssec:composition}).

\subsection{Bilevel Feasibility}\label{ssec:bf}
In this and the next two sections, we only consider the single-commodity case.
Denote \(\setX\) the feasible set of the follower's problem. Note that \(\setX\) is independent from \(t\).
Let \(x\) be a point in \(\setX\). For each \(t \geq 0\), the cost of \(x\) is \(\cx + \txA\).
We call \(g = \cx\) the \emph{base cost}, and \(w = \xA\) the \emph{reduced reaction}. If two points in \(\setX\) share the same base cost
and reduced reaction, they are equivalent in the NPP, because they produce the same follower's cost and leader's revenue regardless of \(t\).
Let \(\setX^* = \{(\cx, \xA) \mid x \in \setX\}\) be the set of all \((g, w)\) pairs that are feasible.
Then, each pair of \((g, w) \in \setX^*\) represents an equivalence class of reactions.

Let \(f(t)\) be the function of the follower's cost, \ie the optimal value of \cref{prog:aggregated} restricted to commodity \(k\).
For \(t \ngeq 0\), assign \(f(t) = -\infty\). It is well-known that \(f\) is a concave function.
For the ease of analysis, we investigate \(-f\) instead of \(f\) due to its convexity.
For \(t \geq 0\), \(-f(t)\) is the maximum of \(-g-\tw\) for all pairs of \((g, w) \in \setX^*\).
Another perspective is to say that the epigraph \(\epi(-f)\) is the intersection of all the closed half-spaces of the form:
\begin{equation}
	H(g, w) = \{(t, z) \mid z \geq -g -\tw \}\label{eq:halfplane}
\end{equation}
where \((g, w) \in \setX^*\). Also denote \(H_=(g, w)\) the hyperplane corresponding to \(H(g, w)\).
By construction, a point \(x \in \setX\) is optimal in the follower's problem for some \(t\) if and only if the hyperplane \(H_=(\cx, \xA)\) supports
\(\epi(-f)\). We call such point \emph{bilevel feasible}.

\begin{definition}\label{def:bfpath}
	A reaction \(x\) is \emph{bilevel feasible} when \(H_=(\cx, \xA)\) supports \(\epi(-f)\).
	Otherwise, it is \emph{bilevel infeasible}.
\end{definition}

\begin{lemma}[Bui~\etal~\cite{bui2022}]\label{lem:bf-equivalence}
	A reaction \(x\) is bilevel feasible if and only if there exists \(t \geq 0\) such that \(x\) is optimal in the follower's problem,
	\ie \(f(t) = \cx + \txA\).
\end{lemma}
\begin{proof}
	\(H_=(\cx, \xA)\) supports \(\epi(-f)\) \(\Leftrightarrow\) There exists \((t, z)\) such that \(z = -\cx -\txA\) and \(z = -f(t)\)
	\(\Leftrightarrow\) \(-f(t) = -\cx -\txA\).
\end{proof}

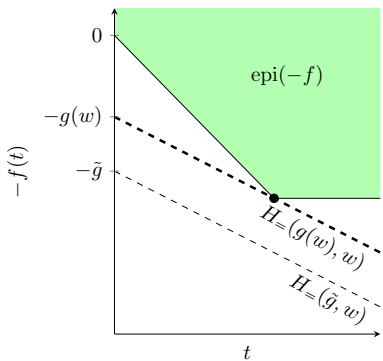
\begin{figure}
	\centering
	\scalebox{0.8}{\begin{tikzpicture}
    \tikzstyle{label}=[inner sep=5pt]
    \tikzstyle{nodelabel}=[inner sep=0]
    \tikzstyle{soliddot}=[only marks,mark=*]
    \tikzstyle{hollowdot}=[fill=white,only marks,mark=*]
    \begin{axis}[
            width=6cm,
            height=7cm,
            axis lines = left,
            ylabel = $-f(t)$,
            xlabel = $t$,
            xmin=0,
            xmax=5,
            ymin=-5.5,
            ymax=0.5,
            xtick=\empty,
            ytick={0,-1.5,-2.5},
            yticklabels={0,$-g(w)$,$-\tilde g$},
            axis on top,
        ]

        \fill[green!30] (0,0) -- (3,-3) -- (5, -3) -- (5, 2) -- (0, 2) -- cycle;
        \draw (0,0) -- (3,-3) -- (5, -3);
        \node at (3.25, -0.75) {$\epi(-f)$};

        \addplot[domain=0:5, dashed, very thick] {-0.5*x-1.5} node[rotate=-28, pos=1.0, below left] {$H_=(g(w), w)$};
        \addplot[domain=0:5, dashed] {-0.5*x-2.5} node[rotate=-28, pos=1.0, below left] {$H_=(\tilde g, w)$};

        \addplot[soliddot] coordinates{(3,-3)};

    \end{axis}
\end{tikzpicture}}
	\caption{\(H_=(g(w), w)\) is the highest hyperplane given \(w\).}
	\label{fig:supporting-hyperplane}
\end{figure}

Given \(w\), in order for a hyperplane \(H_=(g, w)\) to support \(\epi(-f)\), \(g\) must be the smallest possible value such that \((g, w) \in \setX^*\), so that
\(H_=(g, w)\) is the highest of all hyperplanes with slope \(w\).
Thus, we can optimize \(g\) given \(w\). Define
\begin{align*}
	g(w) & = \sup_t \left\{\tw - (-f(-t)) \right\}                                  \\
	     & = \sup_t \left\{f(t) - \tw \right\}. \numberthis{}\label{eq:conjugate-g}
\end{align*}
\Cref{eq:conjugate-g} tells us that \(g(w)\) is the smallest \(g\) such that \(g \geq f(t) - \tw\) for all \(t\), or equivalently,
\(-f(t) \geq -g - \tw\). Recall from \cref{eq:halfplane} that \(H_=(g(w), w)\) supports \(\epi(-f)\) (see \cref{fig:supporting-hyperplane}).
The function \(g(w)\) is called the \emph{convex conjugate} of \(-f(-t)\) \citep{rockafellar1970}.
Since \(-f\) is lower semi-continuous, by the Fenchel-Moreau theorem \citep{rockafellar1970}, the conjugate of the conjugate is the original function, so we also have:
\begin{align*}
	-f(-t)                            & = \sup_w \left\{\tw - g(w) \right\}                                      \\
	\Leftrightarrow \hspace{1cm} f(t) & = \inf_w \left\{\tw + g(w) \right\}. \numberthis{}\label{eq:conjugate-f}
\end{align*}

Comparing \cref{eq:conjugate-g} to \cref{prog:aggregated-dual,prog:conjugate-follower},
we can deduce that \(g(w)\) is the optimal value of \cref{prog:conjugate-follower}.
The same conclusion can be derived from \cref{eq:conjugate-f} and \cref{prog:aggregated,prog:conjugate-dual}.
Therefore, the purpose of the conjugate follower model is to find \(g(w)\), the smallest base cost given a reduced reaction.
This leads us to a redefinition of bilevel feasibility:
\begin{lemma}\label{lem:bf-equivalence2}
	A reaction \(x\) is \emph{bilevel feasible} if and only if \(\cx = g(\xA)\).
\end{lemma}

We consider \cref{lem:bf-equivalence,lem:bf-equivalence2} dual of each other, where one characterizes bilevel feasibility based on \(f(t)\),
and the other based on \(g(w)\). Using the new definition, to test the bilevel feasibility of a particular \(x \in \setX\),
we solve the \cref{prog:conjugate-dual} (or \cref{prog:conjugate-follower})
with \(w\) set to \(\xA\). If the optimal objective value \(g(w)\) is equal to \(\cx\), then it is bilevel feasible.
Otherwise, it is bilevel infeasible.

\begin{example}\label{ex:bf-revisit}
	We revisit \cref{ex:single-comm} with the reaction plot redrawn in \cref{fig:bf-single-reaction}.
	We mentioned the case \(\xA = (0, 1)\) as an example of bilevel infeasible reaction. Indeed, set \(w = (0, 1)\)
	and solve \cref{prog:conjugate-follower}, we get \(g(w) = 5\). The shortest path with \(\xA = (0, 1)\) has a base cost of \(6 > g(w)\),
	thus it is bilevel infeasible.

	In \cref{ex:single-comm}, we only listed 4 binary values of \(\xA\). It raises the question:
	Is there a bilevel feasible \(x\) such that \(\xA\) is non-binary? Solving \cref{prog:conjugate-follower} with \(w = (0.5, 0.5)\)
	produces \(g(w) = 3.5\). In general, if \(w = (\alpha, \alpha)\) where \(0 < \alpha < 1\), then \(g(w) = 5(1 - \alpha) + 2\alpha\).
	\Cref{prog:conjugate-dual} gives us the reaction \(x\) described in \cref{fig:bf-single-value}.
	This is the convex combination of 2 simple paths corresponding to \(\xA = (1, 1)\) and \(\xA = (0, 0)\), respectively.
	By \cref{lem:bf-equivalence}, if this \(x\) is bilevel feasible, then it must be optimal for some \(t\). What is this set of \(t\)?
	\Cref{prog:conjugate-follower} tells us that this reaction is optimal when \(t\) is in the intersection of the
	2 regions labeled as \((1, 1)\) and \((0, 0)\) (the thick segment in \cref{fig:bf-single-reaction}).
	It seems that there is a relation between convex combinations of paths and intersections of regions in the reaction plot,
	which we will explore in the next two sections.

	\begin{figure}
		\centering
		\begin{subfigure}[b]{0.48\textwidth}
			\centering
			\scalebox{0.8}{\begin{tikzpicture}
    \begin{axis}[
            scale=0.6,
            axis equal image,
            axis lines = center,
            xlabel = $t_1$,
            ylabel = $t_2$,
            xmin=0,
            xmax=4,
            ymin=0,
            ymax=4,
        ]
        
        \draw (3,0) -- (1,2) -- (1,4);
        \draw (0,2) -- (1,2);
        \draw[line width = 1mm] (3,0) -- (1,2);
        \node at (1,1) {$(1,1)$};
        \node at (0.5, 3) {$(1,0)$};
        \node at (2.5, 2.5) {$(0,0)$};
    \end{axis}
\end{tikzpicture}}
			\subcaption{Reaction plot}
			\label{fig:bf-single-reaction}
		\end{subfigure}
		\hfill
		\begin{subfigure}[b]{0.48\textwidth}
			\centering
			\scalebox{0.8}{\begin{tikzpicture}[scale=1.5]
    \tikzstyle{vertex}=[circle,draw,minimum size=20pt,inner sep=0pt]
    \tikzstyle{tolled}=[->, >=latex, dashed]
    \tikzstyle{tollfree}=[->, >=latex]

    \node[vertex] (o1) at (1, 1) {$o$};
    \node[vertex] (u) at (1, 0) {$u$};
    \node[vertex] (v) at (2, 0) {$v$};
    \node[vertex] (p) at (3, 0) {$p$};
    \node[vertex] (q) at (4, 0) {$q$};
    \node[vertex] (d1) at (4, 1) {$d$};

    \draw[tollfree] (v) edge node[above]{$\alpha$} (p);
    \draw[tolled] (u) edge node[below]{$\alpha$} (v) (p) edge node[below]{$\alpha$} (q);
    
    \draw[tollfree] (o1) edge node[left] {$\alpha$} (u) (q) edge node[right] {$\alpha$} (d1);
    \draw[tollfree] (o1) edge node[above] {$1 - \alpha$} (d1);


\end{tikzpicture}}
			\subcaption{Value of \(x\) with \(w = (\alpha, \alpha)\)}
			\label{fig:bf-single-value}
		\end{subfigure}
		\caption{Illustrations for \cref{ex:bf-revisit}.}
	\end{figure}
\end{example}

Bilevel feasibility is useful in the preprocessing step to prune away bilevel infeasible paths.
Bilevel infeasible paths, by definition, are not optimal to the follower's problem for any \(t \geq 0\), thus the follower
never considers them regardless of the leader's decision. This process is demonstrated in Bui~\etal~\cite{bui2022},
and we will reuse it in our experiments described in \cref{sec:experiments}.

\subsection{Action-Reaction Duality}\label{ssec:duality}
The relationship between \(f(t)\) and \(g(w)\) is more than just finding the smallest base cost.
Because \(f(t)\) is a polyhedral function, so is \(g(w)\), meaning their epigraphs can be defined as the intersection
of a finite number of closed half-spaces. We will prove that a dimensional relation holds between the sets of optimal \(t\) and \(w\).
This relation is better illustrated with an example.

\begin{example}\label{ex:action-plot}
	We reuse the graph in \cref{fig:multi-comm-series-network}.
	In \cref{ssec:reaction}, we introduced the reaction plot as a tool to describe the optimal solution of \cref{prog:aggregated} given \(t\).
	The reaction plot of the graph in \cref{fig:multi-comm-series-network} is redrawn in \cref{fig:duality-series-reaction}
	Another way to derive the reaction plot is through the projection of \(\epi(-f)\) to the space of \(t\) (see \cref{fig:duality-projection}).
	Each full-dimensional region in the reaction plot corresponds to a facet of \(\epi(-f)\), and its label is the (negative) gradient of \(-f(t)\)
	when \(t\) is in the interior of the region. We can apply the same process to \(\epi(g)\) as well, and the result is called the \emph{action plot},
	shown in \cref{fig:duality-series-action}. The labels in the action plot are \((t_1, \ldots, t_{|\setA_1|})\) which are the gradients of \(g(w)\),
	similar to the labels in the reaction plot.

	There is a correspondence between ``features'' (defined later) of one plot to those of the other.
	Every label representing a full-dimensional region in one plot corresponds to a vertex of the other plot (the red/green vertices and polygons).
	Edges separating two polygons in one plot become edges connecting the corresponding vertices in the other (the blue edges).
	As a general rule, the sum of the dimensions of a pair of features is always \(|\setA_1|\), which is 2 in this example.

	\begin{figure}
		\begin{minipage}[b][][b]{.62\textwidth}
			\begin{subfigure}[b]{0.48\textwidth}
				\centering
				\scalebox{0.8}{\begin{tikzpicture}
    \begin{axis}[
            scale=0.65,
            axis equal image,
            axis lines = center,
            xlabel = $t_1$,
            ylabel = $t_2$,
            xmin=0,
            xmax=8,
            ymin=0,
            ymax=8,
            axis on top,
        ]
        \fill[green!50] (0,0) -- (5,0) -- (5,4) -- (3,6) -- (0,6) -- cycle;

        \draw (5,0) -- (5,4) -- (8,4);
        \draw (3,8) -- (3,6) -- (0,6);
        \draw (3,6) -- (5,4);
        \node at (2.5,3) {$(1,1)$};
        \node at (1.7,7) {$(1,0)$};
        \node at (6.5,2) {$(0,1)$};
        \node at (5.5,6) {$(0,0)$};

        \node[circle,draw,fill=red!70,inner sep=0pt,minimum size=6pt] at (3,6) {};
        \draw[blue, line width = 1mm] (5,0) -- (5,4);
    \end{axis}
\end{tikzpicture}}
				\subcaption{Reaction plot}
				\label{fig:duality-series-reaction}
			\end{subfigure}
			\hfill
			\begin{subfigure}[b]{0.48\textwidth}
				\centering
				\scalebox{0.8}{\begin{tikzpicture}
    \begin{axis}[
            scale=0.6,
            axis equal image,
            axis lines = center,
            xlabel = $w_1$,
            ylabel = $w_2$,
            xmin=0,
            xmax=1.5,
            ymin=0,
            ymax=1.5,
            axis on top,
        ]

        \fill[red!30] (0,0) -- (1,0) -- (1,1) -- cycle;

        \draw (0,1) -- (2,1);
        \draw (1,0) -- (1,2);
        \draw (0,0) -- (1,1);

        \node at (1.25,0.5) {$(0,6)$};
        \node at (0.7,0.3) {$(3,6)$};
        \node at (0.3,0.7) {$(5,4)$};
        \node at (0.5,1.25) {$(5,0)$};
        \node at (1.25,1.25) {$(0,0)$};

        \draw[blue, line width = 1mm] (1,1) -- (0,1);
        \node[circle,draw,fill=green!70!black,inner sep=0pt,minimum size=6pt] at (1,1) {};
    \end{axis}
\end{tikzpicture}}
				\subcaption{Action plot}
				\label{fig:duality-series-action}
			\end{subfigure}
			\addtocounter{figure}{-1}
			\captionof{figure}{Correspondence of features between the action plot and the reaction plot.}
		\end{minipage}%
		\hfill
		\begin{minipage}[b][][b]{.35\textwidth}
			\centering
			\scalebox{0.75}{\tdplotsetmaincoords{110}{-30}
\begin{tikzpicture}[scale=0.35, tdplot_main_coords]
    \def\maxt{10}
    \def\epishift{15}
    \def\maxf{\epishift+2}

    \colorlet{col11}{blue!70!black}
    \colorlet{col10}{green!70!black}
    \colorlet{col01}{orange!80!black}
    \colorlet{col00}{yellow}
    \colorlet{colvert}{white!70!black}

    \newcommand{\projpair}[4]{
        \coordinate (#1) at (#2,#3,\epishift-#4);
        \coordinate (#1P) at (#2,#3,0);
        \coordinate (#1Q) at (#2,#3,\maxf);
    }
    \newcommand{\projdashed}[1]{
        \draw[dashed,opacity=0.6] (#1) -- (#1P);
    }

    \projpair{O}{0}{0}{0};
    \projpair{A}{0}{6}{6};
    \projpair{B}{3}{6}{9};
    \projpair{C}{5}{4}{9};
    \projpair{D}{5}{0}{5};
    \projpair{AI}{0}{\maxt}{6};
    \projpair{BI}{3}{\maxt}{9};
    \projpair{CI}{\maxt}{4}{9};
    \projpair{DI}{\maxt}{0}{5};
    \projpair{EI}{\maxt}{\maxt}{9};

    \draw[dashed, opacity=0.6] (EI) -- (EIQ) -- (AIQ);
    \draw[dashed, opacity=0.6] (EIQ) -- (DIQ);

    \fill[fill opacity=0.9, fill=colvert!30!black] (AI) -- (AIQ) -- (OQ) -- (O) -- (A) -- cycle;
    \fill[fill opacity=0.9, fill=colvert] (DI) -- (DIQ) -- (OQ) -- (O) -- (D) -- cycle;

    \draw[fill opacity=0.9, fill=col11!80!black] (O) -- (A) -- (B) -- (C) -- (D) -- cycle;
    \draw[fill opacity=0.9, fill=col10!30!black] (AI) -- (A) -- (B) -- (BI);
    \draw[fill opacity=0.9, fill=col01!90!black] (CI) -- (C) -- (D) -- (DI);

    \fill[fill opacity=0.9, fill=col00!70!black] (BI) -- (B) -- (C) -- (CI) -- (EI);
    \draw (BI) -- (B) -- (C) -- (CI);

    \foreach \i in {A,B,C,D,AI,BI,CI,DI,EI}{
        \projdashed{\i}
    }

    \draw[fill opacity=0.8, fill=col11] (OP) -- (AP) -- (BP) -- (CP) -- (DP) -- cycle;
    \draw[fill opacity=0.8, fill=col10] (AIP) -- (AP) -- (BP) -- (BIP);
    \draw[fill opacity=0.8, fill=col01] (CIP) -- (CP) -- (DP) -- (DIP);

    \fill[fill opacity=0.8, fill=col00] (BIP) -- (BP) -- (CP) -- (CIP) -- (EIP);
    \draw (BIP) -- (BP) -- (CP) -- (CIP);

    \draw[thick,->] (0,0,0) -- (\maxt,0,0) node[anchor=north east]{$t_1$};
    \draw[thick,->] (0,0,0) -- (0,\maxt,0) node[anchor=north west]{$t_2$};
    \draw[thick,->] (0,0,0) -- (0,0,\maxf) node[anchor=south]{$-f(t)$};

    \node[white!90!black] at ($0.25*(O)+0.25*(B)+0.25*(C)+0.25*(D)$) {$(1, 1)$};
    \node[white!90!black] at ($0.25*(A)+0.25*(B)+0.25*(AI)+0.25*(BI)$) {$(1, 0)$};
    \node[white!90!black] at ($0.25*(C)+0.25*(D)+0.25*(CI)+0.25*(DI)$) {$(0, 1)$};
    \node[white!0!black] at ($0.25*(BI)+0.25*(CI)+0.25*(EI)+0.25*(C)$) {$(0, 0)$};

    \node[white!90!black] at ($0.33*(AI)+0.33*(AIQ)+0.33*(OQ)$) {$(\infty, 0)$};
    \node[white!0!black] at ($0.33*(DI)+0.33*(DIQ)+0.33*(OQ)$) {$(0, \infty)$};

\end{tikzpicture}}
			\captionof{figure}{Reaction plot as a projection of \(\epi(-f)\).}
			\label{fig:duality-projection}
		\end{minipage}
	\end{figure}
\end{example}

As remarked in \cref{ex:action-plot}, the reaction plot is the projection of \(\epi(-f)\).
Thus, the ``features'' (polygons, edges, vertices) in the reaction plot correspond to faces of \(\epi(-f)\).
Conversely, every non-vertical face of \(\epi(-f)\) corresponds to a feature of the reaction plot
with the same dimension. The same applies to faces of \(\epi(g)\). In this section,
we want to establish a bijection between faces of \(\epi(-f)\)
and faces of \(\epi(g)\), along with their dimensional relation.
Hereafter, let \(n = |\setA_1|\) be the dimension of the space of \(t\) (and of \(w\)).

\begin{definition}\label{def:action-set}
	An \emph{action set}\footnote{In \cref{sec:asymmetry}, we call the set of \(t\) that are optimal to \cref{prog:conjugate-follower}
		action set. The term action set here has the same meaning as in \cref{sec:asymmetry}, despite of their different definitions.}
	is the projection of a non-vertical face of \(\epi(-f)\) to the space of \(t\).
	To put it concretely, a set \(T \subseteq \setreal^n\) is an action set when there exists a face \(F\) of \(\epi(-f)\)
	such that:
	\begin{enumerate}[(i)]
		\item The direction \((t, z) = (0, 1)\) is not in \(\aff(F)\) (\(F\) is non-vertical);
		\item \(T = \{t \in \setreal^n \mid (t, z) \in F \text{ for some \(z\)}\}\) (\(T\) is the projection of \(F\)).
	\end{enumerate}
\end{definition}

Condition (i) implies that given \(t \in T\), there is only one value of \(z\) such that \((t, z) \in F\).
Condition (ii) forces \(z\) to be \(-f(t)\). Thus, given \(T\), we can derive
\(F = \{(t, -f(t)) \in \setreal^n \times \setreal \mid t \in T\}\), making the mapping from \(F\) to \(T\) bijective.
A \emph{reaction set}\footnote{This also has the same meaning as in \cref{sec:asymmetry}, but \(\funcR(t)\) is for \(x\), while \(W\) (and later \(\funcW(t)\)) is for \(w\).}
\(W\) is defined similarly with \(\epi(g)\).

\begin{lemma}\label{lem:action-face}
	Given an action set \(T\), if \(T'\) is a face of \(T\), then \(T'\) is an action set.
\end{lemma}
\begin{proof}
	Let \(F = \{(t, -f(t)) \in \setreal^n \times \setreal \mid t \in T\}\) be the face of \(\epi(-f)\) corresponding to \(T\).
	It is easy to verify that \(t \in \ri(T)\) if and only if \((t, -f(t)) \in \ri(F)\),
	where \(\ri(\cdot)\) denotes the relative interior of a set.
	Since \(T'\) is a face of \(T\), it is the intersection of \(T\) and the half-space \(\{t \mid \tr{\alpha}t \leq \beta\}\)
	for some \(\alpha \in \setreal^n\) and \(\beta \in \setreal\).
	Consider \(F' = \{(t,-f(t)) \in F \mid \tr{\alpha}t \leq \beta\}\). Then, \(t \in T'\) if and only if \((t,-f(t)) \in F'\),
	and since \(T' \cap \ri(T) = \varnothing\), we have \(F' \cap \ri(F) = \varnothing\).
	Thus, \(F'\) is a face of \(F\) which in turn is a non-vertical face of \(\epi(-f)\), meaning \(F'\) is a non-vertical face of \(\epi(-f)\).
	Also, \(T'\) is the projection of \(F'\), making it an action set.
\end{proof}

Given \(t \geq 0\), let \(\funcW(t)\) be the set of \(w\) that are optimal to \cref{prog:aggregated}.
Recall from \cref{eq:conjugate-f} that \(w \in \funcW(t)\) if and only if
\begin{equation}\label{eq:conjugate-duality}
	f(t) = \tw + g(w).
\end{equation}
Recall the definition of \(\funcT(w)\) in \cref{ssec:conjugate} as the set of \(t\) that are optimal to \cref{prog:conjugate-follower},
we also obtain the equality above. Therefore, \(w \in \funcW(t)\) if and only if \(t \in \funcT(w)\).

Given \(t \geq 0\), we rewrite \cref{eq:conjugate-f} as:
\[f(t) = \inf_{w, z} \left\{\tw + z \mid z \geq g(w) \right\} = \inf_{(w, z)} \left\{\tw + z \mid (w, z) \in \epi(g) \right\}.\]
By treating \(\epi(g)\) as the feasible set of some optimization problem,
the set of \((w, z)\) that is optimal forms a non-vertical face of \(\epi(g)\). As a result,
\(\funcW(t)\) is a reaction set. Conversely, given a reaction set \(W\), there always exist \(t \geq 0\) such that \(W = \funcW(t)\),
since every face is optimal for some objective function (which we can adjust with \(t\)).
By a similar reasoning, \(\funcT(w)\) is an action set for all \(w \geq 0\).

Given an action set \(T \subseteq \setreal^n\), define $\funcW(T) = \bigcap_{t \in T} \funcW(t).$

\begin{proposition}
	If \(T\) is an action set, then \(\funcW(T)\) is a reaction set.
\end{proposition}
\begin{proof}
	Define \(\funcG(W) = \{(w, g(w)) \in \setreal^n \times \setreal \mid w \in W\}\) as the face of \(\epi(g)\) corresponding to
	the reaction set \(W\). Consider \(G = \bigcap_{t \in T} \funcG(\funcW(t))\). Since all points in \(\funcG(\funcW(t))\)
	are in the form \((w, g(w))\) for some \(w \geq 0\), so are all points in \(G\).
	Given \(w\), all of the following statements are equivalent:
	\begin{align*}
		w \in \funcW(T)
		\Leftrightarrow w \in \funcW(t), \forall t \in T
		\Leftrightarrow  (w,g(w)) \in \funcG(\funcW(t)), \forall t \in T
		\Leftrightarrow (w,g(w)) \in G.
	\end{align*}
	Thus, \(\funcW(T)\) is the projection of \(G\) onto the \(w\)-space. Because \(\funcG(\funcW(t))\) are all non-vertical faces of \(\epi(g)\),
	their intersection \(G\) is either empty or a non-vertical face of \(\epi(g)\).
	Since \(T\) is an action set, there exists \(w\) such that \(T = \funcT(w)\). This means that
	for all \(t \in T = \funcT(w)\), \(w \in \funcW(t)\), hence \(w \in \funcW(T)\) and \(\funcW(T)\) is non-empty.
	Therefore, \(G\) is non-empty and is a non-vertical face of \(\epi(g)\), meaning that its projection \(\funcW(T)\) is a reaction set.
\end{proof}

Given a reaction set \(W\), define \(\funcT(W)\) similarly. Then, \(\funcT\) is the inverse mapping of \(\funcW\).
\begin{proposition}\label{prop:action-identity}
	Given an action set \(T\), then \(\funcT(\funcW(T)) = T\).
\end{proposition}
\begin{proof}
	For all \(t \in T\), \(w \in \funcW(T)\), since \(\funcW(T) \subseteq \funcW(t)\), we have \(w \in \funcW(t)\) implying \(t \in \funcT(w)\).
	Hence, \(t \in \bigcap_{w \in \funcW(T)} \funcT(w) = \funcT(\funcW(T))\), meaning \(T \subseteq \funcT(\funcW(T))\).

	Conversely, since \(T\) is an action set, there exists \(w \geq 0\) such that \(\funcT(w) = T\), hence
	\(w \in \funcW(t)\) for all \(t \in \funcT(w) = T\), implying \(w \in \funcW(T)\). It follows that \(\funcT(\funcW(T)) \subseteq \funcT(w) = T\)
	(by definition of \(\funcT\)). We conclude that \(\funcT(\funcW(T)) =~T\).
\end{proof}

\begin{corollary}\label{cor:action-subset}
	Given two action sets \(T_1\) and \(T_2\), then \(T_1 \subseteq T_2\) if and only if \(\funcW(T_1) \supseteq \funcW(T_2)\).
\end{corollary}
\begin{proof}
	The forward direction is true due to the definition of \(\funcW\). For the reverse direction,
	\(\funcW(T_1) \supseteq \funcW(T_2)\) implies \(\funcT(\funcW(T_1)) \subseteq \funcT(\funcW(T_2))\). By \cref{prop:action-identity},
	\(T_1 \subseteq T_2\).
\end{proof}

Together, \(\funcW\) and \(\funcT\) define a bijection between the action sets and the reaction sets.
Given \(t \geq 0\), the smallest action set containing \(t\) is \(\funcT(\funcW(t))\). The same applies to all \(w \geq 0\).
This mapping can be extended to a bijection between non-vertical faces of \(\epi(-f)\) and those of \(\epi(g)\).
The dimensional property linking them can be expressed~as:
\begin{theorem}\label{thm:duality-dim}
	Given an action set \(T\), it holds $\dim T + \dim \funcW(T) = n.$
\end{theorem}
\begin{proof}
	Since \(T\) is an action set, there exists \(w\) such that \(T = \funcT(w)\), which implies \(w \in \funcW(T)\).
	Let \(S\) be the subspace parallel to \(\aff(T)\) and \(S^\perp\) be the orthogonal subspace of~\(S\).
	Their dimensional relation is $\dim S^\perp = n - \dim S = n - \dim T$.
	Choose a basis of \(S^\perp\) so that \(S^\perp = \spn\{u_1, u_2, \ldots, u_m\}\) where \(m = \dim S^\perp\).
	For each \(i = 1,\ldots,m\), we will prove that there exists \(\delta_i \in \spn\{u_i\}\), \(\delta_i \neq 0\) such that \(w + \delta_i \in \funcW(T)\).

	Suppose that the above statement is false, meaning there exists \(i\) such that for all \(\delta_i \in \spn\{u_i\} \setminus \{0\}\),
	\(w + \delta_i \notin \funcW(T)\). We can build a sequence \(\delta_i^{(j)} \to 0\) in \(\spn\{u_i\} \setminus \{0\}\) such that
	\(w + \delta_i^{(j)} \geq 0\) so that each \(w + \delta_i^{(j)}\) is contained in some reaction set.\footnote{
		If \(w > 0\), then the sequence \(w + \delta_i^{(j)} \geq 0\) is trivial to construct. If \(w_a = 0\) for some \(a \in \setA_1\),
		then \(t_a\) for some \(t \in T\) can get arbitrary large. Since \(\delta_{i} \in S^\perp\), \(\delta_{i,a} = 0\) and \(w_a + \delta_{i,a} = 0\) for all \(i\)
		which do not violate the condition \(w + \delta_i^{(j)} \geq 0\).}
	Since \(\epi(g)\) is a polyhedron, it has a finite number of faces, hence a finite number of reaction sets.\footnote{As suggested,
		this dimensional property does not hold when \(f\) and \(g\) are not polyhedral functions.}
	Thus, there must be a reaction set \(W'\) that contains infinitely many terms of \(w + \delta_i^{(j)}\).
	Reaction sets are closed (they are projections of faces), hence the limit \(w = \lim (w + \delta_i^{(j)})\) is in \(W'\) as well, implying that \(\funcT(w) = T\) contains \(\funcT(W')\).

	Because \(u_i \in S^\perp\), \(\tr{u_i}(t - t') = 0\) for all \(t, t' \in T\). Choose \(t \in T, t' \in \funcT(W')\),
	then for all \(\delta_i \in \spn\{u_i\} \setminus \{0\}\) (also recall \cref{eq:conjugate-duality}):
	\begin{align*}
		\tr{(w + \delta_i)}(t - t')                       & = \tr{w}(t - t')  = f(t) - f(t')                                     & \text{(since \(t, t' \in \funcT(w)\))} \\
		\Leftrightarrow \quad f(t) - \tr{(w + \delta_i)}t & = f(t') - \tr{(w + \delta_i)}t'. \numberthis\label{eq:duality-proof}
	\end{align*}
	Notice that \(t' \in \funcT(W')\) and \(w + \delta_i^{(j)} \in W'\) for some \(j\),
	hence \(t' \in \funcT(w + \delta_i^{(j)})\) and \(f(t') - \tr{(w + \delta_i^{(j)})}t' = g(w + \delta_i^{(j)})\) (by \cref{eq:conjugate-duality}).
	Thus, the left-hand side of \cref{eq:duality-proof} is also equal to \(g(w + \delta_i^{(j)})\) which implies that
	\(w + \delta_i^{(j)} \in \funcW(t)\) for all \(t \in T\). It follows	that \(w + \delta_i^{(j)} \in \funcW(T)\) which is a contradiction.

	To conclude the proof, since there exists \(\delta_i \neq 0\) such that \(w + \delta_i \in \funcW(T)\) in every direction \(u_i \in S^\perp\),
	\(w\) along with all the \(w + \delta_i\) form an affinely independent set in \(\funcW(T)\), thus \(\dim \funcW(T) \geq m\). For all \(u \notin S^\perp\),
	there exists \(t, t' \in T\) such that \(\tr{u}(t - t') \neq 0\) which implies \(w + u \notin \funcW(T)\).
	Therefore, \(\dim \funcW(T) = m = n - \dim T\).
\end{proof}

\begin{example}\label{ex:duality-single-toll}
	We consider the NPP described in \cref{ex:single-toll}, restricted to commodity \(k = 3\).
	The plots of \(-f(t)\) and \(g(w)\) are shown in \cref{fig:duality-followercost,fig:duality-conjugatecost}.
	Each non-vertical face of \(\epi(-f)\) (encoded as \(A, b, C, d\)) has a counterpart in \(\epi(g)\) (encoded as \(a, B, c, D\)) through the mappings
	\(\funcW\) and \(\funcT\). Consider the 1-dimensional face \(b\) of \(\epi(-f)\) (the red segment in \cref{fig:duality-followercost}), which is the intersection of
	\(\epi(-f)\) and the hyperplane corresponding to \((g, w) = (0, 1)\). We can verify that \(g(1) = 0\) (the red vertex in \cref{fig:duality-conjugatecost}).
	Perturbing \(w\) by any small amount changes the optimal face from \(b\) into either \(A\) or \(C\), hence \(B\) is a 0-dimensional face of \(\epi(g)\).
	In contrast, the 0-dimensional face \(C\) of \(\epi(-f)\) (the blue vertex in \cref{fig:duality-followercost}) corresponds to not only \(w = 0.5\), but also
	any \(w\) in its neighborhood. Informally, we can ``wiggle'' the supporting hyperplane (the blue dashed line in \cref{fig:duality-followercost}) while
	keeping \(C\) as the optimal point, which we cannot do with \(b\). As a consequence, \(c\) is a 1-dimensional face of \(\epi(g)\).
	In both cases, the sum of the dimensions is \(|\setA_1| = 1\).

	\begin{figure}
		\centering
		\begin{subfigure}[b]{0.48\textwidth}
			\centering
			\scalebox{0.8}{\begin{tikzpicture}
    \tikzstyle{label}=[inner sep=5pt]
    \tikzstyle{nodelabel}=[inner sep=0]
    \tikzstyle{soliddot}=[only marks,mark=*]
    \tikzstyle{hollowdot}=[fill=white,only marks,mark=*]
    \begin{axis}[
            width=6cm,
            height=7cm,
            axis lines = left,
            ylabel = $-f(t)$,
            xlabel = $t$,
            xtick distance = 2,
            xmin=0,
            xmax=5,
            ymin=-5,
            ymax=1,
            axis on top,
        ]

        \fill[green!30] (0,0) -- (3,-3) -- (5, -3) -- (5, 2) -- (0, 2) -- cycle;

        \draw (0,0) node[above right] {\(A\)}
            --node[above right, red] {\(b\)}
            (3,-3) node[above=3pt, blue] {\(C\)}
            --node[above] {\(d\)} (5, -3);

        \addplot[domain=0:5, dashed, red] {-x} node[label, pos=0.8, below left] {$w=1$};
        \addplot[domain=0:5, dashed, blue] {-0.5*x-1.5} node[pos=0.47, below left] {$w=0.5 \pm \delta$};

        \draw[red, ultra thick] (0,0) -- (3, -3);        
        \addplot[soliddot] coordinates{(0,0)};
        \addplot[blue, soliddot] coordinates{(3,-3)};

        \node at (3.25, -0.5) {$\epi(-f)$};

        \draw[blue, Stealth-Stealth] ([shift=(110:0.5cm)]1,-2) arc[radius=0.5, start angle=110, end angle=197];
    \end{axis}
\end{tikzpicture}}
			\subcaption{Plot of \(-f(t)\)}
			\label{fig:duality-followercost}
		\end{subfigure}
		\hfill
		\begin{subfigure}[b]{0.48\textwidth}
			\centering
			\scalebox{0.8}{\begin{tikzpicture}
    \tikzstyle{label}=[inner sep=5pt]
    \tikzstyle{nodelabel}=[inner sep=0]
    \tikzstyle{soliddot}=[only marks,mark=*]
    \tikzstyle{hollowdot}=[fill=white,only marks,mark=*]
    \begin{axis}[
            width=6cm,
            height=7cm,
            axis lines = left,
            ylabel = $g(w)$,
            xlabel = $w$,
            xtick distance = 1,
            xmin=0,
            xmax=1.7,
            ymin=-2,
            ymax=4,
            axis on top,
        ]

        \fill[green!30] (0,3) -- (1,0) -- (2,0) -- (2,4) -- (0,4) -- cycle;

        \draw (0,3) node[above right] {\(D\)}
            --node[above right, blue] {\(c\)}
            (1,0) node[above=3pt, red] {\(B\)}
            --node[above] {\(a\)} (1.7, 0);

        \addplot[domain=0:2, dashed, red] {-1*x+1} node[pos=0.35, below left] {$t=1 \pm \delta$};
        \addplot[domain=0:2, dashed, blue] {-3*x+3} node[label, pos=0.7, below left] {$t=3$};

        \draw[blue, ultra thick] (0,3) -- (1,0);
        \addplot[soliddot] coordinates{(0,3)};
        \addplot[red, soliddot] coordinates{(1,0)};

        \node at (1, 2.5) {$\epi(g)$};

        \draw[red, Stealth-Stealth] ({0.3+0.15*cos(110)},{0.6+0.45*sin(110)})
            arc[x radius=0.15, y radius=0.45, start angle=110, end angle=197];
    \end{axis}
\end{tikzpicture}}
			\subcaption{Plot of \(g(w)\)}
			\label{fig:duality-conjugatecost}
		\end{subfigure}

		\caption{Plots for \cref{ex:duality-single-toll}.}
	\end{figure}
\end{example}

\subsection{Strong Bilevel Feasibility}\label{ssec:strongbf}
\Cref{ex:bf-revisit} shows that fractional reactions \(x\) can be bilevel feasible.
However, we can write fractional reactions as convex combinations of simple paths (which are the extreme points of \(\setX\)),
hence if a fractional reaction is optimal to the leader's problem,
then there must be an optimal simple path as well. Simple paths always have binary reduced reactions \(w = \xA\).
Moreover, as remarked before, the labels \((w_1, \ldots, w_n)\) of full-dimensional regions in the reaction plot are always binary.
From the discussion in \cref{ssec:duality}, these regions are action sets, and
full-dimensional action sets correspond to 0-dimensional reaction sets (\cref{thm:duality-dim}), which, in turn, correspond
to 0-dimensional faces of \(\epi(g)\). Thus, these labels are the projections of extreme points of \(\epi(g)\).
In this section, we explore the relations between simple paths, binary \(w\), and extreme points of \(\setX\) and \(\epi(g)\).

\begin{definition}\label{def:strongbf}
	A reduced reaction \(w\) is \emph{strongly bilevel feasible} when \(\{w\}\) is a reaction set.
	Otherwise, it is \emph{weakly bilevel feasible}.
\end{definition}
\begin{lemma}\label{lem:strongbf-equivalence}
	Given \(w \geq 0\), the following statements are equivalent:
	\begin{enumerate}[(a)]
		\item \(w\) is \emph{strongly bilevel feasible};
		\item \(\{w\}\) is a reaction set;
		\item \(\funcW(\funcT(w)) = \{w\}\)
		\item \((w,g(w))\) is an extreme point of \(\epi(g)\);
		\item \(\dim \funcT(w) = n\).
	\end{enumerate}
\end{lemma}
\begin{proof}
	Direct results of \cref{def:action-set,thm:duality-dim,def:strongbf}.
\end{proof}

Note that strong bilevel feasibility is a property of \(w\), while bilevel feasibility is a property of \(x\).
It is not incidental that extreme points of \(\setX\) are related to extreme points of \(\epi(g)\).

\begin{proposition}\label{prop:extreme-points}
	If \(w\) is strongly bilevel feasible, there exists an extreme point \(x\) of \(\setX\) such that
	\(x\) is bilevel feasible and satisfies \(\xA = w\).
\end{proposition}
\begin{proof}
	First, we prove that there exists a bilevel feasible \(x\) such that \(\xA = w\).
	We solve for \(g(w)\) using \cref{prog:conjugate-dual} to obtain \(x\). Suppose that \(x_a < w_a\) for some \(a \in \setA_1\). By complementary slackness, all \(t \in \funcT(w)\), which are optimal solutions of \cref{prog:conjugate-follower}, must have \(t_a = 0\). Consequently, \(\dim \funcT(w) < n\). Thus, by \cref{lem:strongbf-equivalence}, \(w\) is weakly bilevel feasible. Therefore, if \(w\) is strongly bilevel feasible, all optimal solutions of \cref{prog:conjugate-dual} must satisfy \(\xA = w\).

	Now that we have \(x\) that is bilevel feasible and \(\xA = w\), we will construct an extreme point \(\hx\) of \(\setX\) that also satisfies these properties.
	Suppose that \(x\) is not an extreme point, we can write \(x\) as a strict combination of extreme points \(x^{(i)}\) and extreme rays.\footnote{
		In the context of the shortest path problem, extreme points correspond to simple paths and extreme rays correspond to loops.}
	Choose any \(i\), let \(\hx = \xA^{(i)}\) and \(\hw = \hxA\).
	For all \(t \in \funcT(w)\), \(t\) makes \(x\) optimal to the follower's problem. Then, \(t\) must make all \(x^{(i)}\) optimal (since the combination is strict),
	hence \(\hx\) is bilevel feasible and \(t \in \funcT(\hw)\). Following the same reasoning in the first part of the proof, we conclude that \(\hxA = \hw = w\).
\end{proof}

As a result, all strongly bilevel feasible \(w\) are binary since extreme points of \(\setX\) are simple paths and they are binary.
Similar to the fact that there is an extreme points of \(\setX\) that is optimal to the follower's problem, there is an
extreme point of \(\epi(g)\) that is optimal to the leader's problem.

\begin{theorem}\label{thm:strongbf-optimal}
	There exists a strongly bilevel feasible \(w\) that is optimal to the NPP, \ie it solves \cref{prog:conjugate-leader}.
\end{theorem}
\begin{proof}
	Let \((w',t)\) be an optimal solution of \cref{prog:conjugate-leader}.
	In order to be an optimal solution, first, \(t \in \funcT(w')\), hence \(w' \in \funcW(t)\).
	If \(w'\) is not an extreme point of \(\funcW(t)\), we write \(w'\) as a strict combination of extreme points \(w^{(i)}\)
	and extreme rays. Since \(w^{(i)} \in \funcW(t)\), the pairs \((w^{(i)}, t)\) are also feasible to \cref{prog:conjugate-leader}.
	Because \(\tr{w'}t\) is the highest revenue to the NPP, so are \(\tr{(w^{(i)})}t\) for all \(i\),
	thus \(w^{(i)}\) are optimal to \cref{prog:conjugate-leader}. Extreme points of \(\funcW(t)\) are 0-dimensional faces of a reaction set,
	hence they are 0-dimensional reaction sets (\cref{lem:action-face}). Therefore, \(w^{(i)}\) are strongly bilevel feasible.
\end{proof}

\Cref{prop:extreme-points,thm:strongbf-optimal} imply that NPP can be solved by enumerating only simple paths.
When we extend this concept to multi-commodity problems in \cref{sssec:composition-strongbf}, strong bilevel feasibility
proves to be a much more powerful property than just a tool to eliminate fractional reactions.

\subsection{Composition}\label{ssec:composition}
In this section, we return to the multi-commodity case and generalize the concepts introduced in the previous sections via
the process of composition. Given a commodity \(k \in \setK\), we add the superscript \(k\) to all notations that are related to this commodity
(including \(\setX^k, \xk, \wk, f^k, g^k, \funcT^k, \funcW^k\)). The objective function \(f\) of the aggregated follower's problem
(\cref{prog:aggregated}) is the sum of all \(f^k\):
\[f(t) = \sumk f^k(t).\]

The conjugate \(g(w)\) is defined as in \cref{eq:conjugate-g}.
Since the conjugate operator of summation is infimal convolution \citep{rockafellar1970}, we also have:
\begin{equation}\label{eq:inf-conv}
	g(w) = \inf\left\{\sumk g^k(w^k) ~\middle|~ \sumk w^k = w \right\}.
\end{equation}

Denote \(\setX = \setX^1 \times \cdots \times \setX^{|\setK|}\) the feasible set of \cref{prog:aggregated}
and \((\xk) = (x^1, \ldots, x^{|\setK|}) \in \setX\) a feasible point of that program. Hereafter,
we call tuples such as \((\xk) = (x^1, \ldots, x^{|\setK|})\) and \((\wk) = (w^1, \ldots, w^{|\setK|})\) \emph{compositions} or \emph{composed reactions}.

\subsubsection{Bilevel Feasibility}\label{sssec:composition-bf}
With the new definitions of \(f(t)\) and \(g(w)\) for the multi-commodity case, \cref{lem:bf-equivalence2} is straightforward to extend:
\begin{definition}\label{def:bf-composed}
	A composed reaction \((\xk) \in \setX\) is \emph{bilevel feasible} when
	\begin{equation}\label{eq:bf-composed}
		\sumk \cxk = g\left(\sumk \xkA\right).
	\end{equation}
\end{definition}

We can still use the conjugate model (\cref{prog:conjugate-follower,prog:conjugate-dual}) to test bilevel feasibility.
The only change is to set \(w = \sumk\xkA\). This will be the common theme of this section.
\Cref{lem:bf-equivalence} is stated similarly. What are more interesting are the conditions that it implies
when we decompose \(f\) and \(g\) into \(f^k\) and~\(g^k\).

\begin{proposition}\label{prop:bf-composed-equivalence}
	Given a composed reaction \((\xk) \in \setX\), the following statements are equivalent:
	\begin{enumerate}[(a)]
		\item \((\xk)\) is bilevel feasible;
		\item \Cref{eq:bf-composed} is satisfied;
		\item Each \(\xk\) is bilevel feasible individually, and $g\left(\sumk \xkA\right) = \sumk g^k(\xkA)$;
		\item There exists \(t \geq 0\) such that \(f(t) = \sumk \left(\cxk + \txkA\right)\);
		\item There exists \(t \geq 0\) such that for all \(k \in \setK\), \(f^k(t) = \cxk + \txkA\).
	\end{enumerate}
\end{proposition}
\begin{proof}
	(a) \(\Leftrightarrow\) (b) is \cref{def:bf-composed}. (a) \(\Leftrightarrow\) (d) is \cref{lem:bf-equivalence} rewritten with the new~\(f\).
	(d) \(\Leftarrow\) (e) is trivial. Since \(f^k(t) \leq \cxk + \txkA\) for all \(k \in \setK\), (d) \(\Rightarrow\)~(e).

	(b) \(\Rightarrow\) (c): Since \(\cxk \geq g^k(\xkA)\), from \cref{eq:inf-conv},
	\begin{equation}\label{eq:bf-composed-proof}
		\sumk \cxk \geq \sumk g^k(\xkA) \geq g\left(\sumk \xkA\right).
	\end{equation}
	If (b) is true, then the inequalities become equalities, and \(\cxk = g^k(\xkA)\) means that
	\(\xk\) is bilevel feasible individually (\cref{lem:bf-equivalence2}). Thus, (c) is true.

	(b) \(\Leftarrow\) (c): Each \(\xk\) is bilevel feasible individually, hence \(\cxk = g^k(\xk)\). Then,
	all inequalities in \cref{eq:bf-composed-proof} are equalities which implies (b).
\end{proof}

Statements (c) and (e) are the decomposed forms of statements (b) and (d), respectively.
Statements (b) and (d) are considered dual of each other, as well as statements (c) and (e).
Statement (e) means that if \((\xk)\) is bilevel feasible, then there exists \(t \geq 0\) such that
all \(\xk\) are simultaneously optimal to their own followers' problems. Statement (c) gives us an extra condition that a bilevel feasible \((\xk)\) needs
to satisfy even when all \(\xk\) are bilevel feasible individually. In what follows, let \(\wk = \xkA\).

\begin{definition}
	A composition \((\wk)\) is \emph{bilevel feasible} when
	\begin{equation}\label{eq:bf-w}
		g\left(\sumk\wk\right) = \sumk g^k(\wk).
	\end{equation}
\end{definition}

Then, statement (c) of \cref{prop:bf-composed-equivalence} can be rewritten as: each \(\xk\) is bilevel feasible individually
and the composition \((\wk)\) is bilevel feasible where \(\wk = \xkA\).
Note that although, by convention, each \(\wk\) is bilevel feasible, their composition may not. In terms of \(\wk\), statement (e) can be restated
exactly the same as statement (c) but with an alternative definition of bilevel feasible~\((\wk)\).

\begin{lemma}\label{lem:bf-action-intersection}
	A composition \((\wk)\) is bilevel feasible if and only if
	\[\funcT\left(\sumk \wk\right) = \bigcap_{k \in \setK} \funcT^k(\wk).\]
\end{lemma}
\begin{proof}
	Given \((\wk)\), for all \(t \geq 0\), we have:
	\begin{equation}\label{eq:bf-w-equivalence}
		\sumk g^k(\wk) \overset{\text{(a)}}\geq g\left(\sumk\wk\right) \overset{\text{(b)}}\geq f(t) - \sumk \twk = \sumk \left(f^k(t) - \twk\right).
	\end{equation}

	(\(\Rightarrow\)) Assuming \cref{eq:bf-w}, then the inequality at (a) is an equality.
	For all \(t \in \funcT\left(\sumk \wk\right)\), the inequality at (b) is satisfied with equality, hence \(g^k(\wk) = f^k(t) - \twk\) and
	\(t \in \funcT^k(\wk)\) for all \(k \in \setK\). For all \(t \in \bigcap_{k \in \setK} \funcT^k(\wk)\),
	all inequalities in \cref{eq:bf-w-equivalence} become equalities and the equality at (b) implies \(t \in \funcT\left(\sumk \wk\right)\).

	(\(\Leftarrow\)) Since \(\funcT\left(\sumk \wk\right)\) is not empty, there exists \(t \geq 0\) such that \(t \in \funcT^k(\wk)\) for all \(k \in \setK\),
	hence \(g^k(\wk) = f^k(t) - \twk\).	Then, all inequalities in \cref{eq:bf-w-equivalence} become equalities and the equality at (a)
	implies \cref{eq:bf-w}.
\end{proof}

The interpretation of \cref{lem:bf-action-intersection} is that a composition \((\wk)\) is bilevel feasible if and only if
all action sets of its components \(\wk\) in the individual reaction plots intersect. Recall from \cref{ex:multi-comm} that the compositions
\(w^1 = (1, 1)\) (\cref{fig:multi-comm-series-reaction}) and \(w^2 = (0, 1)\) (\cref{fig:multi-comm-parallel-reaction})
are bilevel feasible because their action sets intersect to form the aggregated action set with label \(w = (1, 2)\) (\cref{fig:multi-comm-reaction}).
On the other hand, the composition \(w^1 = (0, 1)\) and \(w^2 = (1, 1)\) is bilevel infeasible since their action sets do not intersect despite having the same \(w\).

\subsubsection{Strong Bilevel Feasibility}\label{sssec:composition-strongbf}
All concepts and propositions in \cref{ssec:duality,ssec:strongbf}  are defined based solely on \(f(t)\) and \(g(w)\),
thus they are still valid to the multi-commodity case (except for \cref{prop:extreme-points} where we need to replace \(w = \xA\) with \(w = \sumk\xkA\)).
Given \(w \geq 0\), we call \((\wk)\) a \emph{decomposition} of \(w\) if \((\wk)\) is a composition that satisfies:
\[w = \sumk \wk.\]
A \emph{bilevel feasible decomposition} of \(w\) is a decomposition which is also bilevel feasible.
By \cref{eq:inf-conv}, any \(w \geq 0\) has at least one bilevel feasible decomposition.

\begin{lemma}\label{lem:convex-comb-bf}
	Given \(w \geq 0\), the set of all bilevel feasible decompositions of \(w\) is convex.
\end{lemma}
\begin{proof}
	Let \((\hwk)\) and \((\cwk)\) be two bilevel feasible decompositions of \(w\) and \((\wk)\) be a convex combination of them, \ie
	\(\wk = \lambda\hwk + (1 - \lambda)\cwk\) for all \(k \in \setK\) and some \(0 \leq \lambda \leq 1\). Clearly, \(\sumk \wk = w\),
	hence \((\wk)\) is a decomposition of \(w\). Since \(g^k\) are convex,
	\[\sumk g^k(\wk) \leq \sumk \left(\lambda g^k(\hwk) + (1 - \lambda)g^k(\cwk) \right) = \lambda g(w) + (1 - \lambda) g(w) = g(w).\]
	Because \(\sumk g^k(\wk) \geq g(w)\) (by \cref{eq:inf-conv}), we conclude that \(\sumk g^k(\wk) = g(w)\) and \((\wk)\) is bilevel feasible.
\end{proof}

Similar to the discussion in \cref{sssec:composition-bf}, even when all \(\wk\) are strongly bilevel feasible individually,
it is not guaranteed that \(w\) would be strongly bilevel feasible. Moreover, \(w\) has infinitely many decompositions \((\wk)\),
and not all of them will satisfy the above property. Then, under what conditions is \(w\) strongly bilevel feasible in relation to its
decomposition \((\wk)\)?

\begin{proposition}\label{prop:strongbf-decompose}
	A reduced reaction \(w \geq 0\) is strongly bilevel feasible if and only if
	for all bilevel feasible decompositions \((\wk)\) of \(w\), each \(\wk\) is strongly bilevel feasible individually.
\end{proposition}
\begin{proof}
	(\(\Rightarrow\)) For all bilevel feasible decompositions \((\wk)\) of \(w\), \cref{lem:bf-action-intersection} states that
	\(\funcT(w) \subseteq \funcT^k(\wk)\) for all \(k \in \setK\). If \(w\) is strongly bilevel feasible, then \(\dim \funcT(w) = n\) (\cref{lem:strongbf-equivalence}),
	hence \(\dim \funcT^k(\wk) = n\) for all \(k \in \setK\). Thus, each \(\wk\) is strongly bilevel feasible individually.

	(\(\Leftarrow\)) Suppose that \(w\) is weakly bilevel feasible, we need to show that there exists some bilevel feasible decomposition \((\wk)\)
	such that some \(\wk\) are weakly bilevel feasible individually. By \cref{lem:strongbf-equivalence},
	\((w, g(w))\) is not an extreme point of \(\epi(g)\), hence it can be written as a strict convex combination of
	two distinct points \((\hw, \hat{g})\) and \((\cw, \check{g})\) of \(\epi(g)\):
	\begin{align}
		w    & = \lambda \hw + (1 - \lambda) \cw, \label{eq:strongbf-decompose-1}            \\
		g(w) & = \lambda \hat{g} + (1 - \lambda) \check{g},  \label{eq:strongbf-decompose-2}
	\end{align}
	for some \(0 < \lambda < 1\). Since \(g(w)\) is convex, \(g(w) \leq \lambda g(\hw) + (1 - \lambda)g(\cw)\);
	and because \((\hw, \hat{g})\) and \((\cw, \check{g})\) are in \(\epi(g)\), \(g(\hw) \leq \hat{g}\) and \(g(\cw) \leq \check{g}\).
	It follows that \(g(\hw) = \hat{g}\), \(g(\cw) = \check{g}\), and \(\hw \neq \cw\).
	Let \((\hwk)\) and \((\cwk)\) be some bilevel feasible decompositions of \(\hw\) and \(\cw\), respectively.
	Define \(\wk = \lambda \hwk + (1 - \lambda) \cwk\) for each \(k \in \setK\). \Cref{eq:strongbf-decompose-1}
	implies that \((\wk)\) is a decomposition of \(w\), while \cref{eq:strongbf-decompose-2} implies:
	\begin{align*}
		g(w) & = \lambda g(\hw) + (1 - \lambda) g(\cw) = \sumk \left[ \lambda g^k(\hwk) + (1 - \lambda) g^k(\cwk) \right]                                    \\
		     & \geq \sumk g^k(\wk).                                                                                       &  & \text{(since \(g^k\) convex)}
	\end{align*}
	Because \(\sumk g^k(\wk) \geq g(w)\) (by \cref{eq:inf-conv}), \((\wk)\) is bilevel feasible.
	Moreover, \(\lambda g^k(\hwk) + (1 - \lambda) g^k(\cwk) = g^k(\wk)\) for all \(k \in \setK\).
	Since \(\hw \neq \cw\), there exists \(k \in \setK\) such that \(\hwk \neq \cwk\). Thus,
	\((\wk, g^k(\wk))\) is a strict convex combination of two distinct points \((\hwk, g^k(\hwk))\) and \((\cwk, g^k(\cwk))\)
	of \(\epi(g^k)\). Therefore, \((\wk, g^k(\wk))\) is not an extreme point of \(\epi(g^k)\) and \(\wk\) is weakly bilevel feasible.
\end{proof}

\begin{lemma}\label{lem:uniqueness}
	If \(w \geq 0\) is strongly bilevel feasible, then \(w\) has exactly one bilevel feasible decomposition.
\end{lemma}
\begin{proof}
	Suppose that \(w\) is strongly bilevel feasible and it has more than one bilevel feasible decompositions.
	\Cref{lem:convex-comb-bf} implies that there is an infinite number of bilevel feasible decompositions \((\wk)\) of \(w\),
	hence an infinite number of compositions \((\wk)\) such that each \(\wk\) is strongly bilevel feasible individually
	(\cref{prop:strongbf-decompose}). However, each \(\epi(g^k)\) has a finite number of extreme points, by \cref{lem:strongbf-equivalence},
	the number of strongly bilevel feasible \(\wk\) is finite, so is the number of their compositions \((\wk)\)	which is a contradiction.
\end{proof}

Taking \cref{prop:strongbf-decompose,lem:uniqueness} together, we have:
\begin{theorem}\label{thm:strongbf-uniqueness}
	A reduced reaction \(w \geq 0\) is strongly bilevel feasible if and only if both of the following conditions are satisfied:
	\begin{enumerate}[(i)]
		\item There exists exactly one bilevel feasible decomposition \((\wk)\) of \(w\);
		\item In the unique decomposition \((\wk)\), each \(\wk\) is strongly bilevel feasible individually.
	\end{enumerate}
\end{theorem}

According to \cref{thm:strongbf-uniqueness}, if \(w\) is weakly bilevel feasible and \((\wk)\) is some bilevel feasible decomposition of \(w\),
it might occur that each \(\wk\) is strongly bilevel feasible individually. In this case, however, \(w\) must have another bilevel feasible decomposition.
Such case usually happens when many commodities share the same portion of the graph, hence they indirectly affect each other.
We provide two examples to demonstrate this phenomenon.

\begin{example}\label{ex:degeneracy}
	Consider the multi-commodity NPP with the graph in \cref{fig:degeneracy-network}.
	Commodity 1 can use both tolled arcs in series while commodity 2 can only use the second tolled arc.
	\Cref{fig:degeneracy-reaction,fig:degeneracy-decomposition} illustrate
	the reaction plots of the whole problem and of individual commodities. In \cref{fig:degeneracy-reaction}, notice that there are 2 edges that
	overlap each other (highlighted as a thick line).
	This is due to the structure \(\{u, v, p\}\) being shared by both commodities.
	The label jumps from \((1,0)\) above the edge to \((1,2)\) below the edge,
	skipping the label \((1,1)\). The action plot (\cref{fig:degeneracy-action}) confirms that \(w = (1,1)\) is weakly bilevel feasible,
	whose action set \(\funcT(w)\) is the thick line in \cref{fig:degeneracy-reaction}.
	We can decompose \(w = (1,1)\) in two (over infinitely many) different ways: (i) The red decomposition: \(w^1 = (1, 1)\) and \(w^2 = (0, 0)\); (ii) The blue decomposition: \(w^1 = (1, 0)\) and \(w^2 = (0, 1)\).
	\Cref{fig:degeneracy-decomposition} shows that in both decompositions,
	each \(w^k\) is strongly bilevel feasible individually. However, since there are two bilevel feasible decompositions,
	by \cref{thm:strongbf-uniqueness}, \(w = (1, 1)\) fails to be strongly bilevel feasible.

	The paths corresponding to these decompositions are illustrated in the last two columns of \cref{fig:degeneracy-decomposition}
	(unused arcs are hidden). In both decompositions, both commodities pass through nodes \(u\) and \(v\) but they travel on different subpaths between \(u\) and \(v\):
	one takes the tolled arc \(u-v\), while the other takes the toll-free segment \(u-p-v\). However, a subpath of the shortest path is the shortest subpath.
	Thus, both commodities should use the same subpath between \(u\) and \(v\): either both take the tolled arc, or both take the toll-free segment.
	Dewez~\etal~\cite{dewez2008} use this observation to add an inequality ruling this case out. In our framework developed so far, we state
	that \(w = (1, 1)\) is weakly bilevel feasible, and by \cref{thm:strongbf-optimal}, we should not consider \(w = (1, 1)\)
	as a solution of the NPP in favor of other strongly bilevel feasible values of \(w\) such as \(w = (1, 2)\) (both take the tolled arc) and \(w = (1, 0)\)
	(both take the toll-free segment).

	\begin{figure}
		\centering

		\begin{subfigure}[b]{0.4\textwidth}
			\centering
			\begin{tikzpicture}[scale=1]
    \small
    \tikzstyle{vertex}=[circle,draw,minimum size=16pt,inner sep=0pt]
    \tikzstyle{tolled}=[->, >=latex, dashed]
    \tikzstyle{tollfree}=[->, >=latex]

    \node[vertex] (u) at (-0.25, 0) {$u$};
    \node[vertex] (v) at (1.25, 0) {$v$};
    \node[vertex] (p) at (0.5, -1) {$p$};
    \node[vertex] (o1) at (-1, 1) {$o^1$};
    \node[vertex] (d1) at (2, 1) {$d^1$};
    \node[vertex] (o2) at (-1, -1) {$o^2$};
    \node[vertex] (d2) at (2, -1) {$d^2$};

    \draw[tolled] (o1) edge node[left]{$t_1$} (u);
    \draw[tolled] (u) edge node[above]{$t_2$} (v);
    \draw[tollfree] (o1) edge node[above] {3} (d1);
    \draw[tollfree] (u) edge node[left]{2} (p) (p) edge (v);
    \draw[tollfree] (v) edge (d1);
    \draw[tollfree] (o2) edge (u) (v) edge (d2);

\end{tikzpicture}
			\vspace{0.5cm}
			\subcaption{Graph}
			\label{fig:degeneracy-network}
		\end{subfigure}
		\hfill
		\begin{subfigure}[b]{0.3\textwidth}
			\centering
			\scalebox{0.8}{\begin{tikzpicture}
    \begin{axis}[
            scale=0.6,
            axis equal image,
            axis lines = center,
            xlabel = $t_1$,
            ylabel = $t_2$,
            xmin=0,
            xmax=3.5,
            ymin=0,
            ymax=3.5,
        ]
        \draw (0,2) -- (4,2);
        \draw (1,2) -- (3,0);
        \draw (1,2) -- (1,4);
        \draw[line width = 1mm] (0,2) -- (1,2);

        \node at (1,1) {$(1,2)$};
        \node at (2.7,1.3) {$(0,1)$};
        \node at (0.5,2.7) {$(1,0)$};
        \node at (2.2,2.7) {$(0,0)$};
    \end{axis}
\end{tikzpicture}}
			\subcaption{Overall reaction plot}
			\label{fig:degeneracy-reaction}
		\end{subfigure}
		\hfill
		\begin{subfigure}[b]{0.25\textwidth}
			\centering
			\scalebox{0.8}{\begin{tikzpicture}
    \begin{axis}[
            scale=0.6,
            axis equal image,
            axis lines = center,
            xlabel = $w_1$,
            ylabel = $w_2$,
            xmin=0,
            xmax=1.8,
            ymin=0,
            ymax=2.7,
            xtick distance = 1,
            axis on top,
        ]

        \draw (1,0) -- (1,3);
        \draw (0,1) -- (1,2) -- (3,2);

        \node at (0.5,0.7) {$(1,2)$};
        \node at (1.4,1) {$(0,2)$};
        \node at (0.5,2.1) {$(3,0)$};
        \node at (1.4,2.35) {$(0,0)$};
    \end{axis}
\end{tikzpicture}}
			\subcaption{Action plot}
			\label{fig:degeneracy-action}
		\end{subfigure}
		\par\medskip

		\begin{subfigure}[b]{\textwidth}
			\centering\footnotesize
			\renewcommand{\arraystretch}{1.5}
			\begin{tabular}{r:c:c:c}
				\Xhline{2\arrayrulewidth}
				          & Reaction plots                            & Red decomposition                    & Blue decomposition                    \\
				\hline
				\(k = 1\) & \begin{tikzpicture}[baseline=(current bounding box.center)]
    \begin{axis}[
            scale=0.3,
            axis equal image,
            axis lines = center,
            xmin=0,
            xmax=3.5,
            ymin=0,
            ymax=3.5,
            axis on top,
            every tick label/.append style={font=\scriptsize},
        ]
        \def\labelsize{\scriptsize}

        \fill[red!30] (0,2) -- (1,2) -- (3,0) -- (0,0) -- cycle;
        \fill[blue!30] (0,2) -- (1,2) -- (1,4) -- (0,4) -- cycle;

        \draw (0,2) -- (1,2) -- (3,0);
        \draw (1,2) -- (1,4);

        \node at (1,1) {\labelsize $(1,1)$};
        \node[rotate=90] at (0.5,2.7) {\labelsize $(1,0)$};
        \node at (2.2,2.2) {\labelsize $(0,0)$};
    \end{axis}
\end{tikzpicture} & \begin{tikzpicture}[scale=0.6, baseline=(current bounding box.center)]
    \tikzstyle{vertex}=[circle,draw,minimum size=8pt,inner sep=0pt]
    \tikzstyle{tolled}=[->, >=latex, dashed]
    \tikzstyle{tollfree}=[->, >=latex]
    \tikzstyle{ignored}=[black!30,line width=0.0mm]
    
    \node[vertex, red] (u) at (-0.25, 0) {};
    \node[vertex, red] (v) at (1.25, 0) {};
    \node[vertex] (p) at (0.5, -1) {};
    \node[vertex, red] (o1) at (-1, 1) {};
    \node[vertex, red] (d1) at (2, 1) {};
    \node[vertex] (o2) at (-1, -1) {};
    \node[vertex] (d2) at (2, -1) {};

    \draw[tolled, red] (o1) edge (u);
    \draw[tolled, red] (u) edge (v);
    \draw[tollfree, red] (v) edge (d1);

\end{tikzpicture} & \begin{tikzpicture}[scale=0.6, baseline=(current bounding box.center)]
    \tikzstyle{vertex}=[circle,draw,minimum size=8pt,inner sep=0pt]
    \tikzstyle{tolled}=[->, >=latex, dashed]
    \tikzstyle{tollfree}=[->, >=latex]
    \tikzstyle{ignored}=[black!30,line width=0.0mm]
    
    \node[vertex, blue] (u) at (-0.25, 0) {};
    \node[vertex, blue] (v) at (1.25, 0) {};
    \node[vertex, blue] (p) at (0.5, -1) {};
    \node[vertex, blue] (o1) at (-1, 1) {};
    \node[vertex, blue] (d1) at (2, 1) {};
    \node[vertex] (o2) at (-1, -1) {};
    \node[vertex] (d2) at (2, -1) {};

    \draw[tolled, blue] (o1) edge (u);
    \draw[tollfree, blue] (u) edge (p) (p) edge (v);
    \draw[tollfree, blue] (v) edge (d1);

\end{tikzpicture} \\
				\hdashline
				\(k = 2\) & \begin{tikzpicture}[baseline=(current bounding box.center)]
    \begin{axis}[
            scale=0.3,
            axis equal image,
            axis lines = center,
            xmin=0,
            xmax=3.5,
            ymin=0,
            ymax=3.5,
            axis on top,
            every tick label/.append style={font=\scriptsize},
        ]
        \def\labelsize{\scriptsize}

        \fill[red!30] (0,2) -- (4,2) -- (4,4) -- (0,4) -- cycle;
        \fill[blue!30] (0,2) -- (4,2) -- (4,0) -- (0,0) -- cycle;

        \draw (0,2) -- (4,2);

        \node at (1.7,1) {\labelsize $(0,1)$};
        \node at (1.7,2.7) {\labelsize $(0,0)$};
    \end{axis}
\end{tikzpicture} & \begin{tikzpicture}[scale=0.6, baseline=(current bounding box.center)]
    \tikzstyle{vertex}=[circle,draw,minimum size=8pt,inner sep=0pt]
    \tikzstyle{tolled}=[->, >=latex, dashed]
    \tikzstyle{tollfree}=[->, >=latex]
    \tikzstyle{ignored}=[black!30,line width=0.0mm]
    
    \node[vertex, red] (u) at (-0.25, 0) {};
    \node[vertex, red] (v) at (1.25, 0) {};
    \node[vertex, red] (p) at (0.5, -1) {};
    \node[vertex] (o1) at (-1, 1) {};
    \node[vertex] (d1) at (2, 1) {};
    \node[vertex, red] (o2) at (-1, -1) {};
    \node[vertex, red] (d2) at (2, -1) {};

    \draw[tollfree, red] (u) edge (p) (p) edge (v);
    \draw[tollfree, red] (o2) edge (u) (v) edge (d2);

\end{tikzpicture} & \begin{tikzpicture}[scale=0.6, baseline=(current bounding box.center)]
    \tikzstyle{vertex}=[circle,draw,minimum size=8pt,inner sep=0pt]
    \tikzstyle{tolled}=[->, >=latex, dashed]
    \tikzstyle{tollfree}=[->, >=latex]
    \tikzstyle{ignored}=[black!30,line width=0.0mm]
    
    \node[vertex, blue] (u) at (-0.25, 0) {};
    \node[vertex, blue] (v) at (1.25, 0) {};
    \node[vertex] (p) at (0.5, -1) {};
    \node[vertex] (o1) at (-1, 1) {};
    \node[vertex] (d1) at (2, 1) {};
    \node[vertex, blue] (o2) at (-1, -1) {};
    \node[vertex, blue] (d2) at (2, -1) {};

    \draw[tolled, blue] (u) edge (v);
    \draw[tollfree, blue] (o2) edge (u) (v) edge (d2);

\end{tikzpicture} \\
				\Xhline{2\arrayrulewidth}
			\end{tabular}
			\subcaption{Individual reaction plots and decompositions of \(w = (1, 1)\)}
			\label{fig:degeneracy-decomposition}
		\end{subfigure}

		\caption{Illustrations for \cref{ex:degeneracy}.}
	\end{figure}
\end{example}

\begin{example}\label{ex:degeneracy-triple}
	Strong bilevel feasibility is not limited to the case described in \cref{ex:degeneracy},
	where two commodities passing through a pair of nodes must share the same subpath between those nodes. We can craft
	an even more intricate example involving 3 commodities where an integral weakly bilevel feasible \(w\) only appears when all 3 are considered.
	Consider the NPP with graph in \cref{fig:degeneracy-triple-network}. Commodities 1 and 2 can only
	use tolled arcs 1 and 2 respectively, while commodity 3 can use both in parallel.
	The reaction plots of the overall problem and individual commodities are shown in \cref{fig:degeneracy-triple-reaction,fig:degeneracy-triple-reaction-individual}.
	We are interested in \(w = (1, 1)\) which is weakly bilevel feasible according to the action plot in \cref{fig:degeneracy-triple-action}.
	Its action set \(\funcT(w)\) consists of only a single point \(t = (2, 3)\).
	Similar to \cref{ex:degeneracy}, \(w = (1, 1)\) has two bilevel feasible decompositions \((\wk)\), in which each \(\wk\) is strongly bilevel feasible individually: 	(i) The red decomposition: \(w^1 = (1, 0)\), \(w^2 = (0, 0)\), and \(w^3 = (0, 1)\); (ii) The blue decomposition: \(w^1 = (0, 0)\), \(w^2 = (0, 1)\), and \(w^3 = (1, 0)\).
	The paths of the decompositions are plotted in \cref{fig:degeneracy-triple-decompositions}. In contrast to \cref{ex:degeneracy},
	no pair of nodes is shared between any pair of commodities, which proves that strong bilevel feasibility
	is a more general property than the one described in \cref{ex:degeneracy} (shared pair of nodes).
	Moreover, if any commodity is excluded, this case will disappear, which shows that strong bilevel feasibility involves
	multiple commodities simultaneously and cannot be considered separately.

	\begin{figure}
		\centering

		\begin{subfigure}[b]{0.37\textwidth}
			\centering
			\begin{tikzpicture}[scale=1.2]
    \small
    \tikzstyle{vertex}=[circle,draw,minimum size=16pt,inner sep=0pt]
    \tikzstyle{tolled}=[->, >=latex, dashed]
    \tikzstyle{tollfree}=[->, >=latex]
    \pgfmathsetmacro\yspace{0.6}
    \pgfmathsetmacro\mspace{1.2}
    \pgfmathsetmacro\nspace{\mspace+2*\yspace}

    \node[vertex] (u) at (-\mspace/2, \yspace) {$u$};
    \node[vertex] (v) at (+\mspace/2, \yspace) {$v$};
    \node[vertex] (p) at (-\mspace/2, -\yspace) {$p$};
    \node[vertex] (q) at (+\mspace/2, -\yspace) {$q$};
    \node[vertex] (o1) at (-\nspace/2, 2*\yspace) {$o^1$};
    \node[vertex] (d1) at (+\nspace/2, 2*\yspace) {$d^1$};
    \node[vertex] (o2) at (-\nspace/2, -2*\yspace) {$o^2$};
    \node[vertex] (d2) at (+\nspace/2, -2*\yspace) {$d^2$};
    \node[vertex] (o3) at (-\nspace/2, 0) {$o^3$};
    \node[vertex] (d3) at (+\nspace/2, 0) {$d^3$};

    \draw[tolled] (u) edge node[above]{$1 + t_1$} (v);
    \draw[tolled] (p) edge node[below]{$t_2$} (q);

    \draw[tollfree] (o1) edge (u) (v) edge (d1);
    \draw[tollfree] (o1) edge node[above] {3} (d1);

    \draw[tollfree] (o2) edge (p) (q) edge (d2);
    \draw[tollfree] (o2) edge node[below] {3} (d2);

    \draw[tollfree] (o3) edge (u) (v) edge (d3);
    \draw[tollfree] (o3) edge (p) (q) edge (d3);
    \draw[tollfree] (o3) edge node[above] {4} (d3);

\end{tikzpicture}
			\subcaption{Graph}
			\label{fig:degeneracy-triple-network}
		\end{subfigure}
		\hfill
		\begin{subfigure}[b]{0.3\textwidth}
			\centering
			\scalebox{0.8}{\begin{tikzpicture}
    \tikzstyle{soliddot}=[only marks,mark=*]
    \begin{axis}[
            scale=0.7,
            axis equal image,
            axis lines = center,
            xlabel = $t_1$,
            ylabel = $t_2$,
            xmin=0,
            xmax=6,
            ymin=0,
            ymax=6,
        ]
        \draw (2,0) -- (2,6);
        \draw (0,3) -- (6,3);
        \draw (0,1) -- (3,4) -- (3,6);
        \draw (3,4) -- (6,4);
        
        \addplot[soliddot] coordinates{(2,3)};

        \node at (0.7,2.6) {$(2,1)$};
        \node at (1.2,1) {$(1,2)$};
        \node at (1,4.5) {$(2,0)$};
        \node at (4,1.5) {$(0,2)$};
        \node[rotate=90] at (2.5,4.8) {$(1,0)$};
        \node at (4.3,3.5) {$(0,1)$};
        \node at (4.5,5) {$(0,0)$};
    \end{axis}
\end{tikzpicture}}
			\subcaption{Overall reaction plot}
			\label{fig:degeneracy-triple-reaction}
		\end{subfigure}
		\hfill
		\begin{subfigure}[b]{0.3\textwidth}
			\centering
			\scalebox{0.8}{\begin{tikzpicture}
    \begin{axis}[
            scale=0.7,
            axis equal image,
            axis lines = center,
            xlabel = $w_1$,
            ylabel = $w_2$,
            xmin=0,
            xmax=2.7,
            ymin=0,
            ymax=2.7,
            xtick distance = 1,
            axis on top,
        ]

        \draw (1,0) -- (0,1);
        \draw (2,0) -- (2,1) -- (1,2) -- (0,2);
        \draw (2,1) -- (3,1);
        \draw (1,2) -- (3,2);
        \draw (1,2) -- (1,3);

        \node at (0.3,0.2) {$(3,4)$};
        \node at (1,1) {$(2,3)$};
        \node at (2.4,0.6) {$(0,3)$};
        \node at (0.5,2.2) {$(2,0)$};
        \node at (2.1,1.5) {$(0,1)$};
        \node at (1.8,2.4) {$(0,0)$};
    \end{axis}
\end{tikzpicture}}
			\subcaption{Action plot}
			\label{fig:degeneracy-triple-action}
		\end{subfigure}
		\par\medskip

		\begin{subfigure}[t]{0.58\textwidth}
			\centering\footnotesize
			\begin{tabular}{ccc}
				\makecell{ \(k = 1\) \\ \scalebox{0.95}{\begin{tikzpicture}[
        baseline=(current bounding box.center),
        trim axis left,
        trim axis right,
    ]
    \begin{axis}[
            scale=0.3,
            axis equal image,
            axis lines = center,
            xmin=0,
            xmax=6,
            ymin=0,
            ymax=6,
            axis on top,
            every tick label/.append style={font=\scriptsize},
        ]
        \def\labelsize{\scriptsize}

        \fill[red!30] (2,0) -- (2,6) -- (0,6) -- (0,0) -- cycle;
        \fill[blue!30] (2,0) -- (2,6) -- (6,6) -- (6,0) -- cycle;

        \draw (2,0) -- (2,6);

        \node[rotate=90] at (1,3) {\labelsize $(1,0)$};
        \node at (4,3) {\labelsize $(0,0)$};
    \end{axis}

\end{tikzpicture}}} &
				\makecell{ \(k = 2\) \\ \scalebox{0.95}{\begin{tikzpicture}[
        baseline=(current bounding box.center),
        trim axis left,
        trim axis right,
    ]
    \begin{axis}[
            scale=0.3,
            axis equal image,
            axis lines = center,
            xmin=0,
            xmax=6,
            ymin=0,
            ymax=6,
            axis on top,
            every tick label/.append style={font=\scriptsize},
        ]
        \def\labelsize{\scriptsize}

        \fill[red!30] (0,3) -- (6,3) -- (6,6) -- (0,6) -- cycle;
        \fill[blue!30] (0,3) -- (6,3) -- (6,0) -- (0,0) -- cycle;

        \draw (0,3) -- (6,3);

        \node at (3,1.5) {\labelsize $(0,1)$};
        \node at (3,4.5) {\labelsize $(0,0)$};
    \end{axis}

\end{tikzpicture}}} &
				\makecell{ \(k = 3\) \\ \scalebox{0.95}{\begin{tikzpicture}[
        baseline=(current bounding box.center),
        trim axis left,
        trim axis right,
    ]
    \begin{axis}[
            scale=0.3,
            axis equal image,
            axis lines = center,
            xmin=0,
            xmax=6,
            ymin=0,
            ymax=6,
            axis on top,
            every tick label/.append style={font=\scriptsize},
        ]
        \def\labelsize{\scriptsize}

        \fill[red!30] (0,1) -- (3,4) -- (6,4) -- (6,0) -- (0,0) -- cycle;
        \fill[blue!30] (0,1) -- (3,4) -- (3,6) -- (0,6) -- cycle;

        \draw (0,1) -- (3,4) -- (3,6);
        \draw (3,4) -- (6,4);
        
        \node at (1.5,4.5) {\labelsize $(1,0)$};
        \node at (4,2) {\labelsize $(0,1)$};
        \node at (4.5,5) {\labelsize $(0,0)$};
    \end{axis}

\end{tikzpicture}}} \\
			\end{tabular}
			\subcaption{Individual reaction plots}
			\label{fig:degeneracy-triple-reaction-individual}
		\end{subfigure}
		\hfill
		\begin{subfigure}[t]{0.4\textwidth}
			\centering\footnotesize
			\begin{tabular}{c:c}
				Red                                                     & Blue \\
				\begin{tikzpicture}[scale=0.7, baseline=(current bounding box.center)]
    \tikzstyle{vertex}=[circle,draw,minimum size=8pt,inner sep=0pt]
    \tikzstyle{tolled}=[->, >=latex, dashed]
    \tikzstyle{tollfree}=[->, >=latex]
    \tikzstyle{ignored}=[black!30,line width=0.0mm]
    
    \pgfmathsetmacro\yspace{0.6}
    \pgfmathsetmacro\mspace{1.2}
    \pgfmathsetmacro\nspace{\mspace+2*\yspace}

    \begin{scope}[red]
        \node[vertex] (u) at (-\mspace/2, \yspace) {};
        \node[vertex] (v) at (+\mspace/2, \yspace) {};
        \node[vertex] (p) at (-\mspace/2, -\yspace) {};
        \node[vertex] (q) at (+\mspace/2, -\yspace) {};
        \node[vertex] (o1) at (-\nspace/2, 2*\yspace) {};
        \node[vertex] (d1) at (+\nspace/2, 2*\yspace) {};
        \node[vertex] (o2) at (-\nspace/2, -2*\yspace) {};
        \node[vertex] (d2) at (+\nspace/2, -2*\yspace) {};
        \node[vertex] (o3) at (-\nspace/2, 0) {};
        \node[vertex] (d3) at (+\nspace/2, 0) {};
    
        \draw[tolled] (u) edge (v);
        \draw[tolled] (p) edge (q);
    
        \draw[tollfree] (o1) edge (u) (v) edge (d1);
    
        \draw[tollfree] (o2) edge (d2);
    
        \draw[tollfree] (o3) edge (p) (q) edge (d3);
    \end{scope}

\end{tikzpicture} &
				\begin{tikzpicture}[scale=0.7, baseline=(current bounding box.center)]
    \tikzstyle{vertex}=[circle,draw,minimum size=8pt,inner sep=0pt]
    \tikzstyle{tolled}=[->, >=latex, dashed]
    \tikzstyle{tollfree}=[->, >=latex]
    \tikzstyle{ignored}=[black!30,line width=0.0mm]
    
    \pgfmathsetmacro\yspace{0.6}
    \pgfmathsetmacro\mspace{1.2}
    \pgfmathsetmacro\nspace{\mspace+2*\yspace}

    \begin{scope}[blue]
        \node[vertex] (u) at (-\mspace/2, \yspace) {};
        \node[vertex] (v) at (+\mspace/2, \yspace) {};
        \node[vertex] (p) at (-\mspace/2, -\yspace) {};
        \node[vertex] (q) at (+\mspace/2, -\yspace) {};
        \node[vertex] (o1) at (-\nspace/2, 2*\yspace) {};
        \node[vertex] (d1) at (+\nspace/2, 2*\yspace) {};
        \node[vertex] (o2) at (-\nspace/2, -2*\yspace) {};
        \node[vertex] (d2) at (+\nspace/2, -2*\yspace) {};
        \node[vertex] (o3) at (-\nspace/2, 0) {};
        \node[vertex] (d3) at (+\nspace/2, 0) {};
    
        \draw[tolled] (u) edge (v);
        \draw[tolled] (p) edge (q);
    
        \draw[tollfree] (o1) edge (d1);
    
        \draw[tollfree] (o2) edge (p) (q) edge (d2);
    
        \draw[tollfree] (o3) edge (u) (v) edge (d3);
    \end{scope}

\end{tikzpicture}       \\
			\end{tabular}
			\subcaption{Decompositions of \(w = (1,1)\)}
			\label{fig:degeneracy-triple-decompositions}
		\end{subfigure}

		\caption{Illustrations for \cref{ex:degeneracy-triple}.}
	\end{figure}
\end{example}

\Cref{ex:degeneracy,ex:degeneracy-triple} illustrate the case where a weakly bilevel feasible \(w\) has a decomposition \((\wk)\)
such that each \(\wk\) is strongly bilevel feasible individually, which implies that the decomposition \((\wk)\) is not unique.
It may also happen that a weakly bilevel feasible \(w\) has a unique bilevel feasible decomposition \((\wk)\). In this case,
there must be some \(\wk\) that is weakly bilevel feasible individually. For instance, \(w = (1, 1)\) in \cref{ex:multi-comm}
has a unique decomposition \(w^1 = (0.5, 0.5)\) and \(w^2 = (0.5, 0.5)\), making it weakly bilevel feasible.
In \cref{ssec:strongbfcut}, we will explain how to use \cref{thm:strongbf-uniqueness} to test strong bilevel feasibility of a given \(w\).

We return to the discussion at the start of this section. There are at most \(2^{|\setA_1||\setK|}\) compositions \((\wk)\), among which only some
are bilevel feasible. However, bilevel feasibility alone is not sufficient to reduce this number to \((|\setK|+1)^{|\setA_1|}\),
since for a given \(w\), there can be more than one bilevel feasible decomposition. \Cref{thm:strongbf-uniqueness} tells us that
such \(w\) is weakly bilevel feasible and by \cref{thm:strongbf-optimal}, we can eliminate all decompositions associated with \(w\).
As a result, there are at most \((|\setK|+1)^{|\setA_1|}\) compositions \((\wk)\) that we need to enumerate, since there are at most
\((|\setK|+1)^{|\setA_1|}\) values of \(w\), and each strong bilevel feasible \(w\) has only one bilevel feasible decomposition \((\wk)\).
Therefore, strong bilevel feasibility is a complete description of the asymmetry in the complexity of the NPP.

\section{Numerical Experiments}\label{sec:experiments}
In this section, we demonstrate that the asymmetry in the complexity of the NPP can be exploited
to solve instances with a high number of commodities more efficiently. More specifically, we compare the performance of an MILP solver
under two different MILPs: a control program similar to \cref{prog:aggregated-1level}, and the same program with additional cuts called
\emph{strong bilevel feasibility cuts}. The programs and the procedure to generate the cuts are explained in details in \cref{ssec:strongbfcut}.
Our experimental results are shown in \cref{ssec:results}.

\subsection{Strong Bilevel Feasibility Cut}\label{ssec:strongbfcut}
The NPP can be formulated in many different ways as systematized in Bui~\etal~\cite{bui2022}. In our experiments, we chose the (PASTD) formulation in \citep{bui2022},
which is recalled here for the sake of completeness:
\begin{subequations}
	\label[program]{prog:pastd}
	\begin{align}
		\max_{t,z,y}\  & \sumk \sumpk \eta^k \tpk \zkp \label{eq:pastd-obj}                                                         \\
		\st            & \sumpk \zkp = 1,                                   &  & k \in \setK, \label[progeq]{eq:pastd-flow}         \\
		               & \Ayk \leq c + t,                                   &  & k \in \setK, \label[progeq]{eq:pastd-dual}         \\
		               & \sumpk \ctpk \zkp = \byk,                          &  & k \in \setK, \label[progeq]{eq:pastd-sd}           \\
		               & t \geq 0, \label{eq:pastd-tpos}                                                                            \\
		               & \zkp \geq 0,                                       &  & k \in \setK, \pk \in \setP^k. \label{eq:pastd-end}
	\end{align}
\end{subequations}
This program is similar to \cref{prog:aggregated-1level} except that the arc representation \(\xk_a\) is replaced by the path representation \(\zkp\).
In \cref{prog:pastd}, \(\setP^k\) is the set of bilevel feasible paths, \(\zkp\) is a binary variable which indicates if the path \(\pk \in \setP^k\) is chosen,
while \(\xpk\) is a constant vector, whose element \(\xpka = 1\) if arc \(a \in \setA\) is in path \(\pk\), 0 otherwise. The program also includes
non-unit demand (or volume) \(\eta^k\) for each commodity.
For the details on the linearization of the bilinear terms \(\tr{t}\xpk\zkp\) and the enumeration of paths to produce \(\setP^k\), refer to Bui~\etal~\cite{bui2022}.

The \cref{prog:pastd} is the control program, on which we will add strong bilevel feasibility cuts. The main idea of the cuts is to eliminate compositions of paths \((\xpk)\)
such that \((\xpk)\) is bilevel infeasible or \(w = \sumk \xpkA\) is weakly bilevel feasible, which is undesirable in either case according to \cref{thm:strongbf-optimal}.
To do so, we need a test for strong bilevel feasibility. To form \(\setP^k\), we enumerate bilevel feasible paths \(\xpk\)
which, in the NPP, are already extreme points of \(\setX^k\). Thus, for all \(k \in \setK\), \(\pk \in \setP^k\), \(\wk = \xpkA\) are, in fact, strongly bilevel feasible individually.
By \cref{thm:strongbf-uniqueness}, the aggregated reaction \(w = \sumk \wk\) will be strongly bilevel feasible if and only if \(w\) has no other bilevel feasible
decomposition besides \((\wk)\). Given a composition of paths \((\xpk)\), first we calculate \(w = \sumk \xpkA\), then we solve the following program:
\begin{subequations}
	\label[program]{prog:sbf-test}
	\begin{align}
		\max_{x}\  & \sumk \sum_{a \in \setA_1} (1-\xpka)\xk_a \label[progeq]{eq:sbf-test-obj}                                                   \\
		\st        & \sumk \xkA \leq w,  \label[progeq]{eq:sbf-test-wcap}                                                                        \\
		           & \Axk = \bk,                                                               &  & k \in \setK,                                 \\
		           & \sumk \cxk \leq \sumk \cpk, \label[progeq]{eq:sbf-test-bf}                                                                  \\
		           & \xk \geq 0,                                                               &  & k \in \setK. \label[progeq]{eq:sbf-test-end}
	\end{align}
\end{subequations}
If \cref{prog:sbf-test} has an objective value strictly greater than 0, then there are 3 cases:
\begin{itemize}
	\item Constraint~\cref{eq:sbf-test-bf} is strict, then \(\sumk \cpk\) is not the smallest sum of base costs, hence the composition \((\xpk)\)
	      is bilevel infeasible (\cref{def:bf-composed}).
	\item Constraint~\cref{eq:sbf-test-bf} is an equality but constraint~\cref{eq:sbf-test-wcap} is strict for some \(a \in \setA_1\), then
	      \(t_a = 0\) for all \(t \in \funcT(w)\) (complementary slackness). This means \(\dim \funcT(w) < n\) implying that \(w\) is weakly bilevel feasible.
	\item Both constraints~\cref{eq:sbf-test-bf,eq:sbf-test-wcap} are equalities,
	      then both \((\xpkA)\) and \((\xkA)\) are bilevel feasible decompositions of \(w\).
	      The objective function~\cref{eq:sbf-test-obj} implies that \((\xpkA)\) and \((\xkA)\) are distinct. Thus, \(w\) is weakly bilevel feasible.
\end{itemize}
Therefore, if \cref{prog:sbf-test} has a non-zero objective value, we can add a cut to eliminate the composition \((\xpk)\).

Removing one composition at a time is inefficient and results in a large number of cuts. Besides, the total number of compositions \((\xpk)\) that we need
to test is exponential. Thus, we limit the test to only pairs of paths. Next, we try to group the pairs of paths that will be eliminated into blocks
so we can produce more efficient cuts. The process is described in \cref{alg:strongbfcut-generation}. 
Line 4 calculates the closeness score
\(d_{ij}\) for each pair of commodities, defined as:
\begin{equation}
	d_{ij} = \frac{|\setA^{k_i} \cap \setA^{k_j}|}{\log^2 (|\setP^{k_i}| |\setP^{k_j}|)} \label{eq:clossness-score}
\end{equation}
where \(\setA^{k_i} \subseteq \setA\) is the set of arcs that is used by at least one path in \(\setP^{k_i}\) (same for \(\setA^{k_j}\)).
The idea is that the more arcs two commodities share, the more bilevel infeasible \((\xp, \xq)\) there will be.
The denominator in \cref{eq:clossness-score} prioritizes pairs with low number of paths,
since there are less compositions \((\xp, \xq)\) to test, hence these pairs are faster to process.
Line 3 limits the generation of cuts to only \(N\) pairs. This parameter \(N\) is used to adjust the amount of cuts added to
the model. Lower values of \(N\) mean fewer cuts and less time spent on generating cuts (and more time spent on actually solving the program).
Lines 5-8 conduct the strong bilevel feasibility test between all pairs
of paths in the given pair of commodities to build the matrix \(h_{pq}\).
Line 9 groups all the compositions \((\xp, \xq)\) that will be eliminated, \ie \(h_{pq} = 1\), into blocks. These blocks are the bicliques of a bipartite graph
whose nodes correspond to paths of each commodity \((\setP^{k_i}, \setP^{k_j})\), and edges connect the pairs of paths that will be eliminated.
Thus, within a biclique \((\hat\setP, \hat\setQ), \hat\setP \subseteq \setP^{k_i}, \hat\setQ \subseteq \setP^{k_j}\),
all \((p, q) \in \hat\setP \times \hat\setQ\) can be removed. Therefore, in line 10, for each biclique \((\hat\setP, \hat\setQ)\), we add
an inequality constraint in the form of:
\begin{equation}
	\sum_{p \in \hat\setP} z^{k_i}_p + \sum_{q \in \hat\setQ} z^{k_j}_q \leq 1. \label{eq:sbf-cut}
\end{equation}
Ideally, we want to add as few cuts as possible, so in line 9, our aim is to find the cover with the fewest number of bicliques.
Unfortunately, the minimum biclique cover problem is NP-hard \citep{fleischner2009}, so we use a heuristic.
Note that the process in \cref{alg:strongbfcut-generation} is just an ad-hoc procedure,
which is used in our experiments to solve some practical issues and only serves as a proof of concept.

\begin{algorithm}
	\caption{Generation of strong bilevel feasibility cuts}
	\label{alg:strongbfcut-generation}

	\begin{algorithmic}[1] \footnotesize
		\Require{The NPP instance, parameter \(N\).}
		\Ensure{\Cref{prog:pastd} with strong bilevel feasibility cuts.}
		\State For each \(k \in \setK\), enumerate the set of bilevel feasible paths \(\setP^k\)
		\State For each \((k_i, k_j) \in \setK^2\), calculate the closeness score \(d_{ij}\)
		\State Choose \(N\) pairs in \(\setK^2\) with the highest \(d_{ij}\) to form \(\setK^2_N\)
		\ForAll{\((k_i, k_j) \in \setK^2_N\)}
		\ForAll{\((p, q) \in \setP^{k_i} \times \setP^{k_j}\)}
		\State Run \cref{prog:sbf-test} on composition \((\xp, \xq)\)
		\State \(h_{pq} \gets 1\) if the objective value is non-zero, 0 otherwise
		\EndFor
		\State Find a solution of the biclique cover problem between \(\setP^{k_i}\) and \(\setP^{k_j}\) with \(h_{pq}\) being
		the adjacency matrix
		\State Add cuts of the form in \cref{eq:sbf-cut} based on the biclique cover
		\EndFor
	\end{algorithmic}
\end{algorithm}

\subsection{Setup and Results}\label{ssec:results}
We tested the performance of the cuts by varying the parameter \(N\) in \cref{alg:strongbfcut-generation} within the set
\(\{0, 10, 20, 50, 100, 200, 500, 1000\}\). The case with no cuts (\(N = 0\)) is the control case. If the cuts are effective,
then we should see better results when \(N\) increases. All cases were tested on 350 randomly generated instances,
whose graphs are \(L \times L\) grid with \(L\) decreasing from 12 to 6 (50 instances for each \(L\)).
Decreasing \(L\) makes the graph smaller, hence lowers the difficulty. To compensate for this, we increase the number of commodities \(|\setK|\)
so that the product \(|\setA_1||K|\) stays constant (which is the number of binary variables of the MILP, as explained in \cref{ssec:reaction}).
The particular values of \(L\) and \(|\setK|\) are listed in \cref{tab:kl-relation}. In the last 5 columns of \cref{tab:kl-relation},
10 instances are generated for each cell. There are 5 columns and 7 rows, which results in 350 instances in total.
Other parameters such as the distribution of the costs, the proportion of tolled arcs, etc., can be found in Bui~\etal~\cite{bui2022}.

\begin{table}
	\caption{The number of commodities \(|\setK|\) in relation to the grid size \(L\).}
	\label{tab:kl-relation}
	\centering \footnotesize
	\begin{tabular}{cccc:ccccc}
		\Xhline{2\arrayrulewidth}
		\(L\) & \(|\setV|\) & \(|\setA|\) & Ratio & \multicolumn{5}{c}{\(|\setK|\)} \\
		\hline
		12    & 144         & 528         & 1.00  & 30  & 35  & 40  & 45  & 50      \\
		11    & 121         & 440         & 1.20  & 36  & 42  & 48  & 54  & 60      \\
		10    & 100         & 360         & 1.47  & 44  & 51  & 59  & 66  & 73      \\
		9     & 81          & 288         & 1.83  & 55  & 64  & 73  & 82  & 92      \\
		8     & 64          & 224         & 2.36  & 71  & 82  & 94  & 106 & 118     \\
		7     & 49          & 168         & 3.14  & 94  & 110 & 126 & 141 & 157     \\
		6     & 36          & 120         & 4.40  & 132 & 154 & 176 & 198 & 220     \\
		\Xhline{2\arrayrulewidth}
	\end{tabular}
\end{table}

We conducted 2800 runs in total (8 values of \(N\) times 350 instances).
At the start of each run, commodity-wise preprocessing (described in \citep{bui2022}) is applied.
Then, strong bilevel feasibility cuts are generated and added to \cref{prog:pastd} based on \(N\).
The preprocessing and cut generation are implemented in Julia.
After all cuts are added, we use Gurobi 1.6.1 \citep{gurobi} to solve the final MILP.
Every step is done in a single thread. The time limit for each run is 1 hour, including the time for preprocessing, cut generation, and solution.
All runs were done on the Narval cluster provided by Compute Canada.

\begin{figure}
	\centering
	\scalebox{0.8}{\input{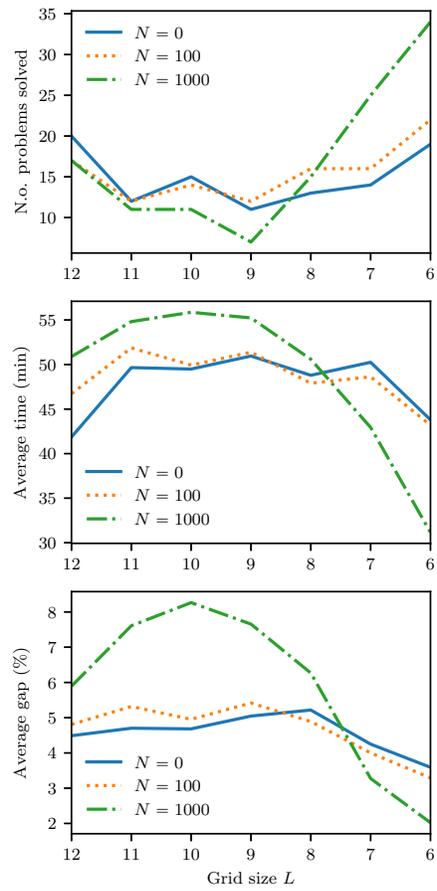}}
	\caption{Plots of the results for selected values of \(N\).}
	\label{fig:results}
\end{figure}

The runs are grouped by \(N\) and \(L\), with 50 runs per group (10 instances for 5 different values of \(|\setK|\) given \(L\)).
We evaluate each group using 3 metrics: the number of instances solved to optimality within the time limit, the average time, and the average optimality gap.
For the first metric, higher is better, while it is the reverse for the other two.
The results are shown in \cref{tab:results}. \Cref{fig:results} plots the metrics for \(N = \) 0, 100, and 1000.
According to the results, increasing \(N\) improves the performance significantly when \(L\) is low. Specifically, at \(L = 6\), the case \(N = 1000\) can solve
15/50 more instances compared to the control case \(N = 0\). On the contrary, the cuts have a detrimental effect on larger graphs, \ie when \(L\) is high.
Thus, the cuts are more effective in the cases with high \(|\setK|\) and low \(|\setA_1|\).
This demonstrates that there exists an asymmetry in the complexity of the NPP and it is possible to take advantage of this asymmetry
to accelerate the solution of the NPP with a small graph and a high number of commodities.

\begin{table}
	\caption{Results of the experiment.}
	\label{tab:results}
	\centering \footnotesize
	\begin{tabular}{cc:ccccccc}
		\Xhline{2\arrayrulewidth}
		\multirow{2}{*}{Metrics} & \multirow{2}{*}{$N$} & \multicolumn{7}{c}{$L$} \\
		\cline{3-9}
		& & 12 & 11 & 10 & 9 & 8 & 7 & 6 \\
		\hline
		\multirow{8}{*}{\makecell{Number of \\ solved \\ instances}}
		& 0 & 20 & 12 & 15 & 11 & 13 & 14 & 19 \\
		& 10 & 19 & 11 & 13 & 11 & 14 & 15 & 19 \\
		& 20 & 19 & 11 & 15 & 11 & 14 & 14 & 21 \\
		& 50 & 19 & 13 & 13 & 11 & 13 & 16 & 20 \\
		& 100 & 17 & 12 & 14 & 12 & 16 & 16 & 22 \\
		& 200 & 17 & 12 & 14 & 10 & 14 & 21 & 25 \\
		& 500 & 17 & 11 & 12 & 9 & 16 & 27 & 31 \\
		& 1000 & 17 & 11 & 11 & 7 & 15 & 25 & 34 \\
		\hline
		\multirow{8}{*}{\makecell{Average \\ time \\ (minutes)}}
		& 0 & 42 & 50 & 49 & 51 & 49 & 50 & 44 \\
		& 10 & 42 & 50 & 49 & 52 & 49 & 50 & 44 \\
		& 20 & 43 & 51 & 50 & 51 & 49 & 50 & 45 \\
		& 50 & 44 & 51 & 50 & 51 & 49 & 50 & 43 \\
		& 100 & 47 & 52 & 50 & 51 & 48 & 49 & 43 \\
		& 200 & 49 & 53 & 52 & 52 & 48 & 48 & 40 \\
		& 500 & 51 & 55 & 54 & 54 & 48 & 42 & 36 \\
		& 1000 & 51 & 55 & 56 & 55 & 51 & 43 & 31 \\
		\hline
		\multirow{8}{*}{\makecell{Average \\ gap \\ (\%)}}
		& 0 & 4.5 & 4.7 & 4.7 & 5.0 & 5.2 & 4.2 & 3.6 \\
		& 10 & 4.5 & 4.9 & 4.6 & 5.3 & 5.1 & 4.2 & 3.6 \\
		& 20 & 4.8 & 4.8 & 4.6 & 4.9 & 5.1 & 4.3 & 3.4 \\
		& 50 & 4.6 & 4.9 & 4.7 & 5.1 & 5.0 & 4.2 & 3.5 \\
		& 100 & 4.8 & 5.3 & 5.0 & 5.4 & 4.9 & 4.0 & 3.3 \\
		& 200 & 5.2 & 5.7 & 5.5 & 5.9 & 5.0 & 3.7 & 3.1 \\
		& 500 & 5.7 & 6.6 & 6.2 & 6.4 & 5.0 & 3.1 & 2.4 \\
		& 1000 & 5.9 & 7.6 & 8.3 & 7.7 & 6.3 & 3.3 & 2.0 \\
		\Xhline{2\arrayrulewidth}
	\end{tabular}
\end{table}

\section{Conclusion}\label{sec:conclusion}
In this paper, we proved that the number of commodities \(|\setK|\) and the number of tolled arcs \(|\setA_1|\) scale the NPP differently.
Particularly, if \(|\setA_1|\) is fixed, the problem can be solved in polynomial time with respect to \(|\setK|\).
We explained the logic behind this asymmetry with the use of convex conjugate, through which we derived the conjugate model,
bilevel feasibility, and strong bilevel feasibility. To the best of our knowledge, this approach is innovative compared to
other applications of the convex conjugate in bilevel programming such as Fenchel-Lagrange duality \citep{aboussoror2017}.
Indeed, these tools can be extended to other pricing problems as long as the assumption that \(f\) is a polyhedral function is kept
(as utilized in \cref{thm:duality-dim}). Future directions include the extension to the case where the follower solves an MILP.

\section*{Acknowledgments}
This work was funded by FRQ-IVADO Research Chair in Data Science for Combinatorial Game Theory, and the NSERC grant 2019-04557.

This research was enabled in part by support provided by Calcul Qu\'ebec (\url{www.calculquebec.ca}) and the Digital Research Alliance of Canada (\url{https://alliancecan.ca/}).

\bibliographystyle{plainnat}
\bibliography{ref}

\begin{thebibliography}{17}
\providecommand{\natexlab}[1]{#1}
\providecommand{\url}[1]{\texttt{#1}}
\expandafter\ifx\csname urlstyle\endcsname\relax
  \providecommand{\doi}[1]{doi: #1}\else
  \providecommand{\doi}{doi: \begingroup \urlstyle{rm}\Url}\fi

\bibitem[Aboussoror et~al.(2017)Aboussoror, Adly, and Saissi]{aboussoror2017}
A.~Aboussoror, S.~Adly, and F.~E. Saissi.
\newblock An extended fenchel--lagrange duality approach and optimality
  conditions for strong bilevel programming problems.
\newblock \emph{SIAM Journal on Optimization}, 27\penalty0 (2):\penalty0
  1230--1255, 2017.
\newblock \doi{10.1137/16M1080896}.

\bibitem[Bouhtou et~al.(2007)Bouhtou, van Hoesel, van~der Kraaij, and
  Lutton]{bouhtou2007}
Mustapha Bouhtou, Stan van Hoesel, Anton~F. van~der Kraaij, and Jean-Luc
  Lutton.
\newblock Tariff optimization in networks.
\newblock \emph{INFORMS Journal on Computing}, 19\penalty0 (3):\penalty0
  458--469, 2007.
\newblock \doi{10.1287/ijoc.1060.0177}.

\bibitem[Brotcorne et~al.(2001)Brotcorne, Labb\'{e}, Marcotte, and
  Savard]{brotcorne2001}
Luce Brotcorne, Martine Labb\'{e}, Patrice Marcotte, and Gilles Savard.
\newblock A bilevel model for toll optimization on a multicommodity
  transportation network.
\newblock \emph{Transportation Science}, 35\penalty0 (4):\penalty0 345--358,
  2001.
\newblock \doi{10.1287/trsc.35.4.345.10433}.

\bibitem[Brotcorne et~al.(2011)Brotcorne, Cirinei, Marcotte, and
  Savard]{brotcorne2011}
Luce Brotcorne, Fabien Cirinei, Patrice Marcotte, and Gilles Savard.
\newblock An exact algorithm for the network pricing problem.
\newblock \emph{Discrete Optimization}, 8\penalty0 (2):\penalty0 246--258,
  2011.
\newblock \doi{10.1016/j.disopt.2010.09.003}.

\bibitem[Brotcorne et~al.(2012)Brotcorne, Cirinei, Marcotte, and
  Savard]{brotcorne2012}
Luce Brotcorne, F~Cirinei, Patrice Marcotte, and Gilles Savard.
\newblock A tabu search algorithm for the network pricing problem.
\newblock \emph{Computers \& Operations Research}, 39\penalty0 (11):\penalty0
  2603--2611, 2012.
\newblock \doi{10.1016/j.cor.2012.01.005}.

\bibitem[Bui et~al.(2022)Bui, Gendron, and Carvalho]{bui2022}
Quang~Minh Bui, Bernard Gendron, and Margarida Carvalho.
\newblock A catalog of formulations for the network pricing problem.
\newblock \emph{INFORMS Journal on Computing}, 34\penalty0 (5):\penalty0
  2658--2674, 2022.
\newblock \doi{10.1287/ijoc.2022.1198}.

\bibitem[Dewez et~al.(2008)Dewez, Labb\'{e}, Marcotte, and Savard]{dewez2008}
Sophie Dewez, Martine Labb\'{e}, Patrice Marcotte, and Gilles Savard.
\newblock New formulations and valid inequalities for a bilevel pricing
  problem.
\newblock \emph{Operations Research Letters}, 36\penalty0 (2):\penalty0
  141--149, 2008.
\newblock ISSN 0167-6377.
\newblock \doi{10.1016/j.orl.2007.03.005}.

\bibitem[Didi-Biha et~al.(2006)Didi-Biha, Marcotte, and Savard]{didibiha2006}
Mohamed Didi-Biha, Patrice Marcotte, and Gilles Savard.
\newblock Path-based formulations of a bilevel toll setting problem.
\newblock In \emph{Optimization with Multivalued Mappings: Theory,
  Applications, and Algorithms}, pages 29--50. Springer US, 2006.
\newblock \doi{10.1007/0-387-34221-4\_2}.

\bibitem[Fleischner et~al.(2009)Fleischner, Mujuni, Paulusma, and
  Szeider]{fleischner2009}
Herbert Fleischner, Egbert Mujuni, Daniël Paulusma, and Stefan Szeider.
\newblock Covering graphs with few complete bipartite subgraphs.
\newblock \emph{Theoretical Computer Science}, 410\penalty0 (21):\penalty0
  2045--2053, 2009.
\newblock ISSN 0304-3975.
\newblock \doi{10.1016/j.tcs.2008.12.059}.

\bibitem[{Gurobi Optimization, LLC}(2022)]{gurobi}
{Gurobi Optimization, LLC}.
\newblock {Gurobi Optimizer Reference Manual}, 2022.
\newblock URL \url{https://www.gurobi.com}.

\bibitem[Karmarkar(1984)]{karmarkar1984}
N.~Karmarkar.
\newblock A new polynomial-time algorithm for linear programming.
\newblock In \emph{Proceedings of the Sixteenth Annual ACM Symposium on Theory
  of Computing}, STOC '84, pages 302--311, New York, NY, USA, 1984. Association
  for Computing Machinery.
\newblock ISBN 0897911334.
\newblock \doi{10.1145/800057.808695}.

\bibitem[Labb{\'e} and Marcotte(2021)]{labbe2021}
Martine Labb{\'e} and Patrice Marcotte.
\newblock \emph{Bilevel Network Design}, pages 255--281.
\newblock Springer International Publishing, Cham, 2021.
\newblock ISBN 978-3-030-64018-7.
\newblock \doi{10.1007/978-3-030-64018-7\_9}.

\bibitem[Labb\'{e} et~al.(1998)Labb\'{e}, Marcotte, and Savard]{labbe1998}
Martine Labb\'{e}, Patrice Marcotte, and Gilles Savard.
\newblock A bilevel model of taxation and its application to optimal highway
  pricing.
\newblock \emph{Management Science}, 44\penalty0 (12-part-1):\penalty0
  1608--1622, 1998.
\newblock \doi{10.1287/mnsc.44.12.1608}.

\bibitem[McCormick(1976)]{mccormick1976}
Garth~P McCormick.
\newblock Computability of global solutions to factorable nonconvex programs:
  Part i—convex underestimating problems.
\newblock \emph{Mathematical programming}, 10\penalty0 (1):\penalty0 147--175,
  1976.
\newblock \doi{10.1007/BF01580665}.

\bibitem[Roch et~al.(2005)Roch, Savard, and Marcotte]{roch2005}
Sébastien Roch, Gilles Savard, and Patrice Marcotte.
\newblock An approximation algorithm for stackelberg network pricing.
\newblock \emph{Networks}, 46\penalty0 (1):\penalty0 57--67, 2005.
\newblock \doi{10.1002/net.20074}.

\bibitem[Rockafellar(1970)]{rockafellar1970}
R.T. Rockafellar.
\newblock \emph{Convex Analysis}.
\newblock Princeton Landmarks in Mathematics and Physics. Princeton University
  Press, 1970.
\newblock ISBN 9780691015866.

\bibitem[{van Hoesel} et~al.(2003){van Hoesel}, {van der Kraaij}, Mannino,
  Oriolo, and Bouhtou]{vanhoesel2003}
C.P.M. {van Hoesel}, A.F. {van der Kraaij}, C.~Mannino, G.~Oriolo, and
  M.~Bouhtou.
\newblock Polynomial cases of the tarification problem.
\newblock Technical report, METEOR, Maastricht University School of Business
  and Economics, January 2003.

\end{thebibliography}

\end{document}